\documentclass[11pt]{article}
\usepackage{amsmath,amsfonts,latexsym}
\usepackage{amscd,amssymb,amsmath}
\usepackage{graphicx}
\usepackage[english]{babel}
\usepackage{indentfirst,enumerate}
\usepackage[T1]{fontenc}
\usepackage[cp1250]{inputenc}
 \usepackage[colorlinks=false]{hyperref}
\usepackage{anysize}
\marginsize{3cm}{2cm}{1cm}{1cm}
\newtheorem{theorem}{Theorem}[section]
\newcounter{tmp}

\newtheorem{lemma}[theorem]{Lemma}

\newtheorem{corollary}[theorem]{Corollary}
\newtheorem{proposition}[theorem]{Proposition}
\newtheorem{definition}[theorem]{Definition}
\newtheorem{remark}[theorem]{Remark}
\newcommand{\hbx}{\hfill$\Box$}
\newcommand{\rmod}{\mathrm{mod}}
\newcommand{\rmi}{\mathrm{i}}

\begin{document}

\title{Displacement sequence of an orientation preserving circle homeomorphism}
\date{}
\author{Wac\l aw Marzantowicz and Justyna Signerska}
\maketitle
\medskip

\begin{abstract}

We give a complete description of the behaviour of the sequence of displacements
$\eta_n(z)=\Phi^n(x) - \Phi^{n-1}(x) \ \rmod \ 1$, $z=\exp(2\pi \rmi x)$,
along a trajectory $\{\varphi^{n}(z)\}$, where $\varphi$ is  an orientation
preserving circle homeomorphism and $\Phi:\mathbb{R} \to \mathbb{R}$ its lift. If the rotation number
$\varrho(\varphi)=\frac{p}{q}$ is rational then $\eta_n(z)$ is
asymptotically  periodic with semi-period $q$. This convergence
to a periodic sequence is uniform in  $z$ if we admit that some points are iterated backward instead of taking only forward iterations for all $z$. If $\varrho(\varphi) \notin \mathbb{Q}$ then the
values of $\eta_n(z)$ are dense in a set which depends on
the  map $\gamma$ (semi-)conjugating $\varphi$ with the rotation
by $\varrho(\varphi)$ and which is the support of the displacements distribution. We provide an effective formula for the displacement distribution if $\varphi$ is $C^1$-diffeomorphism and show approximation of the displacement distribution by sample displacements measured along a trajectory of any other circle homeomorphism which is sufficiently close to the initial homeomorphism $\varphi$. Finally, we prove that even for the irrational rotation number $\varrho$ the displacement sequence exhibits some regularity properties.

\end{abstract}


\section{Introduction}

So far dynamical system theory has been oriented mainly towards studying the distribution of orbits rather than the distribution of displacements along the orbits but the last might be sometimes also a notion of importance. A particular example when the displacement sequence of a circle map is considered, are the so-called \emph{interspike-intervals} for periodically driven \emph{integrate-and-fire} models of neuron's activity (see, for example, \cite{brette1,wmjs1}). In these usually one-dimensional models a continuous dynamics induced by the differential equation is
interrupted by the threshold and reset behaviour,
\begin{eqnarray}
  \qquad \ \ \ \dot{x} &=&f(t, x),  \quad f:\mathbb{R}^2\to\mathbb{R}, \nonumber\\
  \lim_{t\to s^+}x(t) &=& x_r  \qquad \textrm{if} \ x(s)=x_{\Theta}, \nonumber
\end{eqnarray}
 meaning that once a dynamical variable $x(t)$ starting at time $t_0$ from a resting value $x=x_r$ reaches a certain threshold $x_{\Theta}$ at some time $t_1$, it is immediately
reset to a resting value and the system evolves again from a new initial
condition $(x_r,t_1)$ until some time $t_2$ when the threshold is reached again, etc. The question is to describe the sequence of consecutive resets $t_n$ as
iterations of some map $\Phi^n(t_0)$, called the \emph{firing map}, and the sequence of interspike-intervals $t_n-t_{n-1}$ (time intervals between the resets) as a sequence of displacements $\Phi^n(t_0)-\Phi^{n-1}(t_0)$
along a trajectory of this map. The problem appears in various applications,
such as modelling of an action potential (spiking) by a neuron, cardiac rhythms and arrhythmias or electric
discharges in electrical circuits  (see \cite{car-hop} and references therein). Analysis of the behaviour of the displacement sequence of trajectories of an orientation preserving homeomorphism of the circle, covers an answer to this question for the firing map induced by a function $f$ regular enough and periodic in $t$-variable. This special type systems were our motivation for the study presented in this paper. However, our results are of more general character.

We start with homeomorphisms with rational rotation number, where in particular we show the connection between semi-periodic circle homeomorphism and the notion of a semi-periodic sequence and introduce the concept of an $\varepsilon$-basins-shred, separating the points, for which the displacement sequence becomes periodic with given $\varepsilon$-accuracy faster (in terms of number of iterates) when iterated forward than when iterated backward, and the points with the opposite property.

In next we consider homeomorphisms with irrational rotation number. We provide the formula for the displacement distribution with respect to the invariant measure and discuss how this distribution depends on the homeomorphism/diffeomorphism $\varphi$ and the homeomorphism $\gamma$ conjugating it with the rotation $r_{\varrho}$. Finally, with the use of topological dynamics, we show how the recurrent properties of points iterated under $\varphi$ are reflected in the displacement sequence.

Let $\varphi: S^1 \to S^1$ be a map and $\Phi: \mathbb{R} \to \mathbb{R} $ its lift, where
$\mathbb{R}$ covers $S^1$ by the covering projection $\mathfrak{p}:
x\mapsto \exp (2\pi \rmi x)$. If $\varphi:S^1\to S^1$ is an orientation preserving homeomorphism,
then $\Phi(x+1)=\Phi(x)+1$ for all $x\in\mathbb{R}$.

\begin{definition}\label{rotation number}\rm
For $x\in \mathbb{R}$ the limit
\begin{equation}\label{rot nr}
\varrho(\Phi)(x) := \lim_{n\to\infty} \, \frac{\Phi^n(x)}{n}
\end{equation}
is called the \emph{rotation number} of
$\Phi$ at $x$ provided the limit exists.
\end{definition}

\begin{remark}\label{another rotation number}
Let $\Psi(x):=\Phi(x)-x$ be the displacement function associated with $\Phi$. Then
\begin{equation}\label{rotation by displacement}
\varrho(\Phi)(x)= \lim_{N\to\infty}\frac{1}{N}\sum_{n=1}^{N}\Psi(\Phi^{n-1}(x))
\end{equation}
\end{remark}
If $\Phi$ is a lift of an orientation preserving homeomorphism $\varphi\colon S^1\to S^1$, then $\varrho(\Phi)(x)$ exists and does not depend on $x$, following the classical Poincar\'e theory. In this case we define $\varrho(\varphi)\colon=\varrho(\Phi)  \ \rmod  \ 1$, where $\Phi$ is any lift of $\varphi$. Since throughout the rest of the paper we will consider only orientation preserving circle homeomorphisms, we skip the assumption that $\varphi$ preserves orientation in formulation of the forthcoming theorems and definitions.
\begin{definition}\label{ciag disp}
The sequence
\begin{equation}\label{dispseq}
\eta_{n}(z):=\Psi(\Phi^{n-1}(x))\ \rmod \ 1 = \Phi^n(x)-\Phi^{n-1}(x)  \ \rmod \ 1, \quad n=1,2,\ldots
\end{equation} will be called the \emph{displacement sequence} of a point $z=\exp(2\pi \rmi x)\in S^1$.
\end{definition}

Note that $\eta_n(z)$ can be seen as an arc length from the point $\varphi^{n-1}(z)$ to $\varphi^{n}(z)$ with respect to the positive orientation of $S^1$. In particular it does not depend on a choice of the lift $\Phi$.

At first we make two simple observations
 \begin{remark}\label{drobna uwaga}
 If $\varphi$ is a rotation by $2\pi\varrho$, where $\varrho$ can be either rational or irrational, then the sequence $\eta_{n}(z)$ is constant. Precisely, $\eta_n(z)=\varrho$  for  all $z\in S^1$ and $n\in\mathbb{N}$.
 \end{remark}
\begin{remark}\label{drobna uwaga 2}
If $\varphi$ is conjugated to the rational rotation by $2\pi\varrho$, where $\varrho=\frac{p}{q}$, then  $\varphi$ is $q$-periodic, i.e. $\Phi^q(x)=x+p$. Consequently, the sequence $\eta_n(z)$ is $q$-periodic and has the same elements for all $z$.
\end{remark}

\section{Semi-periodic circle homeomorphism}
\subsection{General properties}
We recall definitions of a semi-periodic circle homeomorphism  (after \cite{craciun}) and a semi-periodic sequence  (after \cite{berg}):
\begin{definition} A circle homeomorphism with rational rotation number which is not conjugated to a rotation is called \emph{semi-periodic}.
\end{definition}
\begin{definition}\label{semiperiodic}
A sequence $\{x_n\}$ is \emph{semi-periodic} if
\begin{equation}\label{rnie1}
\forall_{\varepsilon>0}\ \exists_{r\in\mathbb{N}}\ \forall_{n\in\mathbb{N}}\ \forall_{k\in\mathbb{N}} \quad |x_{n+rk}-x_n|<\varepsilon
\end{equation}
\end{definition}
Since we want to investigate asymptotic behaviour of orbits we introduce additionally the concept of asymptotic semi-periodicity:
\begin{definition}\label{asympsemiperiodic}
A sequence $\{x_n\}$ is \emph{asymptotically semi-periodic} if
\begin{equation}\label{rnie2}
\forall_{\varepsilon>0}\ \exists_{N\in\mathbb{N}}\ \exists_{r\in\mathbb{N}}\ \forall_{n>N}\ \forall_{k\in\mathbb{N}}\quad |x_{n+rk}-x_n|<\varepsilon
\end{equation}
\end{definition}
There is also a simpler notion of asymptotic periodicity:
\begin{definition}
We say that a sequence $\{x_n\}$ is \emph{asymptotically periodic} if there exists a periodic sequence $\{a_n\}$ such that $\lim_{n\to\infty}\vert x_n-a_n\vert \to 0$.
\end{definition}
Note that the definition of asymptotic semi-periodicity is more general since for an asymptotically semi-periodic sequence this ``semi-period'' $r$ might depend on $\varepsilon$, whereas it does not for asymptotically periodic one.

In this section we will see that the displacement sequence of a semi-periodic circle homeomorphism is asymptotically periodic, which is a natural consequence of the fact that each non-periodic orbit is attracted to some periodic orbit. Moreover, as we  show in Theorem  \ref{uniwersalne N_1}, for a semi-periodic homeomorphism $\varphi$ with $\varrho(\varphi)=\frac{p}{q}$ and given $\varepsilon>0$, there exists a natural number $N$ such that every point $z\in S^1$ starting from $Nq$-iteration forward or from $Nq$-iteration backward is placed within $\varepsilon$-neighbourhood of a periodic orbit. An analogous property obviously holds for displacement sequences of points under $\varphi$, which is formulated in Proposition \ref{wlasnosc dla przemieszczenia}.
\begin{proposition}\label{conj1}
For a semi-periodic circle homeomorphism $\varphi$ the sequence
$\eta_{n}(z)$ is asymptotically periodic (and thus in particular asymptotically semi-periodic) for any $z\in S^1$. Precisely, if $\varrho(\varphi)=p/q$ then for every $z\in S^1$:
\begin{equation}\label{semiokresowosc}
\forall_{\varepsilon>0}\ \exists_{N\in\mathbb{N}} \ \forall_{n>N}\ \forall_{k\in\mathbb{N}}\ \quad |\eta_{n+kq}(z)-\eta_{n}(z)|<\varepsilon
\end{equation}
\end{proposition}
\noindent{\bf Proof.} For all periodic points the statement reduces to Remark \ref{drobna uwaga 2}. Given a non-periodic point $z=\exp(2\pi\rmi x)\in S^1$ there exists a periodic point $z_0=\exp(2\pi\rmi x_0)\in S^1$ and some $\widetilde{N}$ such that for all $n\geq \widetilde{N}$ and $i=0,1,...,q-1$ we have $\vert\Phi^{nq+i}(x)-\Phi^{nq+i}(x_0)\vert<\varepsilon/4$, i.e. the non-periodic orbit of $z$ is asymptotic to the periodic orbit of $z_0$. Then the property (\ref{semiokresowosc}) of the displacement sequence $\eta_n(z)$  holds for $N:=\widetilde{Nq}$. If $p$ and $q$ are relatively prime, then the ``semi-period'' $r=q$ is minimal.\hbx

 \subsection{A uniform choice of $N(z)$}
Given a homeomorphism $\varphi$ with $\varrho(\varphi)=\frac{p}{q}$ and $\varepsilon>0$, it does not exist $N$ such that for $n>N$ and $k\in\mathbb{N}$ we have $\vert\eta_{n+kq}(z)-\eta_{n}\vert<\varepsilon$ for all $z\in S^1$. Nevertheless, it is possible to find one $N$ that
would fit all the points if we allow that for some points we
consider positive iterates and for some negative.

Suppose that $\varrho(\varphi)=\frac{p}{q}$ (we admit also $q=1$, where periodic points of $\varphi$ are precisely fixed points). If $z_*^-$ and $z_{*}^+$ are consecutive periodic points, then every point
$z\in (z_*^-,z_*^+)$ is forward asymptotic under $\varphi^q$ to $z_*^+$ and backward,
i.e. under $\varphi^{-q}$, asymptotic to $z_*^-$ or the other way around.

Suppose now that the set of periodic points $\textrm{Per}(\varphi)$ is finite and ordered as
$\textrm{Per}(\varphi)=\{z^1,z^2, \dots, z^r\}$, $r\in\mathbb{N}$. Let us fix
$\varepsilon >0$ and $k \in \{1, 2, \dots , r\}$. Assume without the
loss of generality that the two consecutive periodic points $z^k$ and $z^{k+1}$ (which belong to different periodic orbits if there is more than one periodic orbit) are, respectively,
backward and forward attracting under $\varphi^q$ for $z\in (z^k,z^{k+1})$. To clarify the notation  we denote $z^k$ by $z_0^-$, and $z^{k+1}$ by
$z_0^+$. It follows that all the points within the interval $(z_{i}^-,z_{i}^+)$, where $z_{i}^-=\varphi^{i}(z_{0}^-)$, $z_{i}^+=\varphi^{i}(z_{0}^+)$ for $i=0,1,\dots, q-1$ and $z_{q}^-=z_{0}^-, \ z_{q}^+=z_{0}^+$, go forward under $\varphi^q$ to $z_{i}^+$ and backward to $z_{i}^-$.

For a given $m\in \mathbb{N}$ we define the functions $\tau_m^+,\; \tau_m^-:\ [z_0^-,z_0^+] \to [0, \max\limits_{0\leq i \leq q-1}\vert z_{i}^+-z_{i}^-\vert]$:
\[
 \tau_m^+(z):= \max\limits_{0\leq i \leq q-1}\vert \varphi^{mq+i}(z)-z_{i}^+\vert, \quad \tau_m^-(z)\colon=\max\limits_{0\leq i \leq q-1}\vert \varphi^{-mq-i}(z)-z_{q-i}^-\vert
\]
 with the following properties

\begin{tabular}{p{1cm}p{14cm}}

  (i) & $\tau_m^+(z)$ is strictly decreasing, $\tau_m^+(z_0^-)=\max\limits_{0\leq i \leq q-1}\vert z_{i}^+-z_{i}^-\vert$, $\tau_m^+(z_0^+)=0$; \\
  (ii) & $\tau_m^-(z)$ is strictly increasing, $\tau_m^-(z_0^-)= 0 $, $\tau_m^-(z_0^+)=\max\limits_{0\leq i \leq q-1}\vert z_{i}^+-z_{i}^-\vert$; \\
  (iii) & if $m^\prime > m$ then $\tau_{m^\prime}^+(z) < \tau_m^+(z)$ and $\tau_{m^\prime}^-(z) < \tau_m^+(z)$ for every $z\in  (z_0^-,z_0^+)$. \\
  \end{tabular}

For $\varepsilon>0$ and $m\in\mathbb{N}$ denote the subsets of $[z_0^-,z_0^+]$:
$$U^+_m(\varepsilon):= \{z: \, \tau_m^+(z) < \varepsilon \}, \quad U^-_m(\varepsilon):= \{z: \, \tau_m^-(z) < \varepsilon\}$$

Then $z <z^\prime$ and
$z\in U^+_m(\varepsilon)$ implies $z^\prime\in U^+_m(\varepsilon)$ and,
analogously, if $z^\prime < z  $ and $z\in U^-_m(\varepsilon)$, then
$z^\prime \in U^+_m(\varepsilon)$. Put $a_m(\varepsilon):= \inf\{z\in(z_0^{-},z_0^+): \ z\in U^+_m(\varepsilon)\}$  and $b_m(\varepsilon):=\sup\{z\in(z_0^{-},z_0^+): \ z\in U^-_m(\varepsilon)\}$. It is clear that $U_m^+(\varepsilon) = (a_m(\varepsilon), z_0^+]$, $U_m^-(\varepsilon) = [z_0^-, b_m(\varepsilon))$, $\bigcup_{m=1}^{\infty} \overline{U}^+_m(\varepsilon)
= [z_0^-, z_0^+]$ and $\bigcup_{m=1}^{\infty}
\overline{U}_{m(\varepsilon)}  = [z_0^-, z_0^+]$.

For fixed $\varepsilon>0$ there exists $m$ such that
$U_m(\varepsilon)\colon = U^-_m(\varepsilon) \cap  U^+_m(\varepsilon) \neq \emptyset$. Let $\tilde{m}=  \tilde{m}(\varepsilon):= \min\{ m\in \mathbb{N}:
\,U^-_m(\varepsilon) \cap U^+_m(\varepsilon) \neq \emptyset\}.$ Then $${\overline{U}}_{\tilde{m}(\varepsilon)}^+
\cap \; {\overline{U}}_{\tilde{m}(\varepsilon)}^{\,-} =
[a_{\tilde{m}(\varepsilon)}, b_{\tilde{m}(\varepsilon)}]$$

 is a closed interval with nonempty interior. We easily justify:

 \begin{tabular}{p{0.5cm}p{13.5cm}}

  (iv) & For every $z \in  (a_{\tilde{m}(\varepsilon)}, b_{\tilde{m}(\varepsilon)})$
and  $m\geq \tilde{m}(\varepsilon) $ we have $\vert \varphi^{mq+i}(z) - \varphi^{i}(z_0^+)\vert < \varepsilon$ and $\vert \varphi^{-mq-i}(z) - \varphi^{q-i}(z_0^-)\vert  <\varepsilon, $ $i=0,1,\dots, q-1$.\\
  (v) & For every $z\in [z_0^-,a_{\tilde{m}(\varepsilon)})$ and  $m\geq
\tilde{m}(\varepsilon) $ we have $\vert \varphi^{-mq-i}(z) - \varphi^{q-i}(z_0^-)\vert <\varepsilon$, $i=0,1,\dots, q-1$. \\
  (vi) & For every $ z \in (b_{\tilde{m}(\varepsilon)}, z_0^+]$ and  $m\geq
\tilde{m}(\varepsilon) $ we have $\vert \varphi^{mq+i}(z) - \varphi^{i}(z_0^+)\vert<\varepsilon$, $i=0,1,\dots, q-1$. \\
  \end{tabular}

\begin{proposition}\label{separating set}
Let $\varphi: S^1 \to S^1$ be a circle homeomorphism which has finitely many periodic points $\{z^1, z^2, \dots, z^r\}$. Fix $\varepsilon>0$ and consider the interval $(z^k,z^{k+1})$ between the two consecutive
different periodic points $z^k$ and $z^{k+1}$ of $\varphi$.

Suppose that $z^{k+1}$ is attracting (under $\varphi^q$) and $z^k$ is repelling within the interval $(z^k,z^{k+1})$. Then there exists a point $\widetilde{z}_k\in (z^k,z^{k+1})$ with the following properties:

\begin{tabular}{p{0.5cm}p{13.5cm}}

  1) & if $z\in B_k^+:=[\widetilde{z}_k, z^{k+1})$ and for some $n\in \mathbb{N}$ $\vert\varphi^{-nq-i}(z)-\varphi^{q-i}(z^k)\vert<\varepsilon$ for all $i=0,1,\dots, q-1$, then also $\vert\varphi^{nq+i}(z)-\varphi^i(z^{k+1})\vert<\varepsilon$, $i=0,1\dots, q-1$; \\
  2) & if $z\in B_k^-:=(z^k,\widetilde{z}_k]$ and for some $n\in\mathbb{N}$ $\vert\varphi^{nq+i}(z)-\varphi^i(z^{k+1})\vert<\varepsilon$ for all $i=0,1,\dots, q-1$, then also $\vert\varphi^{-nq-i}(z)-\varphi^{q-i}(z^k)\vert<\varepsilon$, $i=0,1\dots, q-1$.\\
  \end{tabular}

If the point $z^k$ is attracting and $z^{k+1}$ is repelling, then $B_k^+:=(z^k,\widetilde{z}_k]$,  $B_k^-:=[\widetilde{z}_k, z^{k+1})$ and the analogues of $1)$ and $2)$ hold.

The same occurs if there is only one periodic, i.e. fixed, point $z_0$ but then
$$S^1\setminus \{z_0\}= B_0^+ \cup B_0^- \;\;{\textrm{and}}\;\;\; B_0^{\;+}\cap B_0^{\;-} = \{\tilde{z}_0\}\,.$$
\end{proposition}
The above proposition says that for every point $z\in B_k^+$ its positive semi-orbit $\{\varphi^n(z)\}_{n\in\mathbb{N}}$ in shorter time (in terms of number of iterates) is placed in the $\varepsilon$-neighbourhood of the periodic orbit $\{z^{k+1}, \varphi(z^{k+1}), \dots, \varphi^{q-1}(z^{k+1})\}$ (i.e. for sufficiently large $n$ $\vert \varphi^{nq+i}(z)-\varphi^i(z)\vert<\varepsilon$ for every $i=0,1,\dots, q-1$) than its negative semi-orbit $\{\varphi^{-n}(z)\}_{n\in\mathbb{N}}$ is placed in the $\varepsilon$-neighbourhood of the repelling orbit $\{z^k, \varphi^{-1}(z^k), \dots, \varphi^{-(q-1)}(z^k)\}$. Similarly, the orbits of points of $B_k^-$ are faster, but in negative time, attracted to the $\varepsilon$-neighbourhood of the orbit of $z^k$ than  to the $\varepsilon$-neighbourhood of the orbit of $z^{k+1}$.

\begin{definition}
We call a one-point set $\{\tilde{z}_k\}$  the {\it$\varepsilon$-basins-shred},
since it divides the whole basin $B_k$ into the positive and
negative sub-basins $B_k^+$ and $B_k^-$, respectively, and is a
common border of them.
\end{definition}

\noindent{\bf{Proof of Proposition \ref{separating set}.}} Let us
consider the interval $[a_{\tilde{m}}, b_{\tilde{m}}]\subset (z^k, z^{k+1})$. By
definition $\tau_{\tilde{m}}^+:  [a_{\tilde{m}}, b_{\tilde{m}}] \to
[0,\varepsilon]$ with $\tau_{\tilde{m}}^+(a_{\tilde{m}})= \varepsilon$.
Correspondingly, $\tau_{\tilde{m}}^-: [a_{\tilde{m}}, b_{\tilde{m}}] \to
[0,\varepsilon]$ with $\tau_{\tilde{m}}^-(b_{\tilde{m}})= \varepsilon$. Moreover, $\tau_{\tilde{m}}^-(a_{\tilde{m}}) <  \tau_{\tilde{m}}^-(b_{\tilde{m}})=\tau_{\tilde{m}}^+(a_{\tilde{m}})$ and $\tau_{\tilde{m}}^+(b_{\tilde{m}}) <
\tau_{\tilde{m}}^+(a_{\tilde{m}})=\tau_{\tilde{m}}^-(b_{\tilde{m}}).$ There exists a unique point $\tilde{z}_k \in
(a_{\tilde{m}}, b_{\tilde{m}}) $ such that $ \tau_{\tilde{m}}^+(\tilde{z}_k)
= \tau_{\tilde{m}}^-(\tilde{z}_k).$ Now the statement of Proposition \ref{separating set} follows from
the properties $(iv)-(vi)$. \hbx

\begin{remark}\rm
Note that for a given $\varepsilon$-accuracy  the basins-shred
$\tilde{z}_k \in (z_k, z_{k+1})$ is defined by $\tau_{\tilde{m}_k}^+$ and
$\tau_{\tilde{m}_k}^-$ and  thus it depends on $\tilde{m}_k$. Since
$\tilde{m}_k$ depends on  $\varepsilon$, the basins-shred
$\tilde{z}_k$ changes if we change the $\varepsilon$-accuracy of
approximation.
\end{remark}

\begin{figure}
  \includegraphics[width=0.56\textwidth]{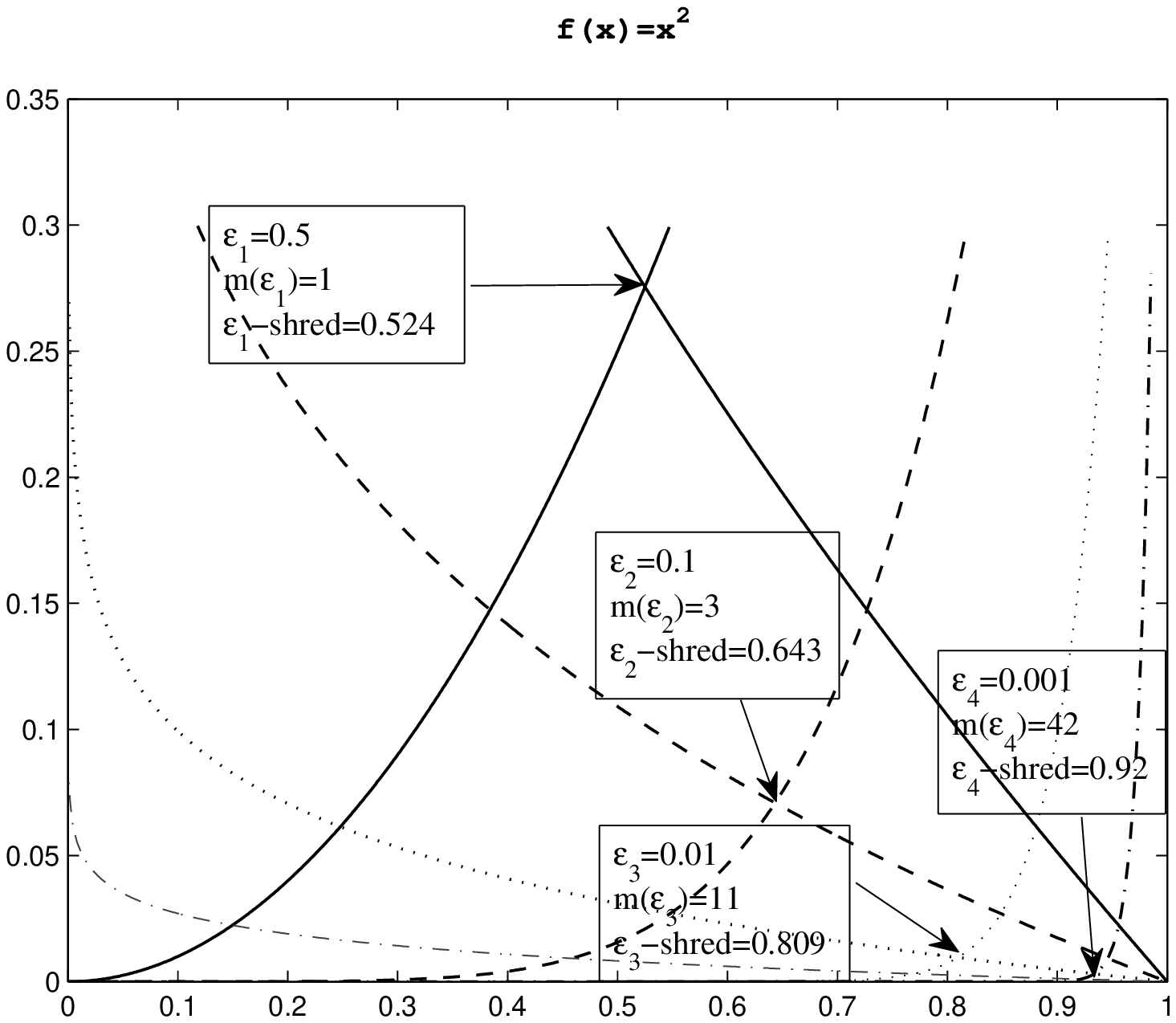}
  \includegraphics[width=0.56\textwidth]{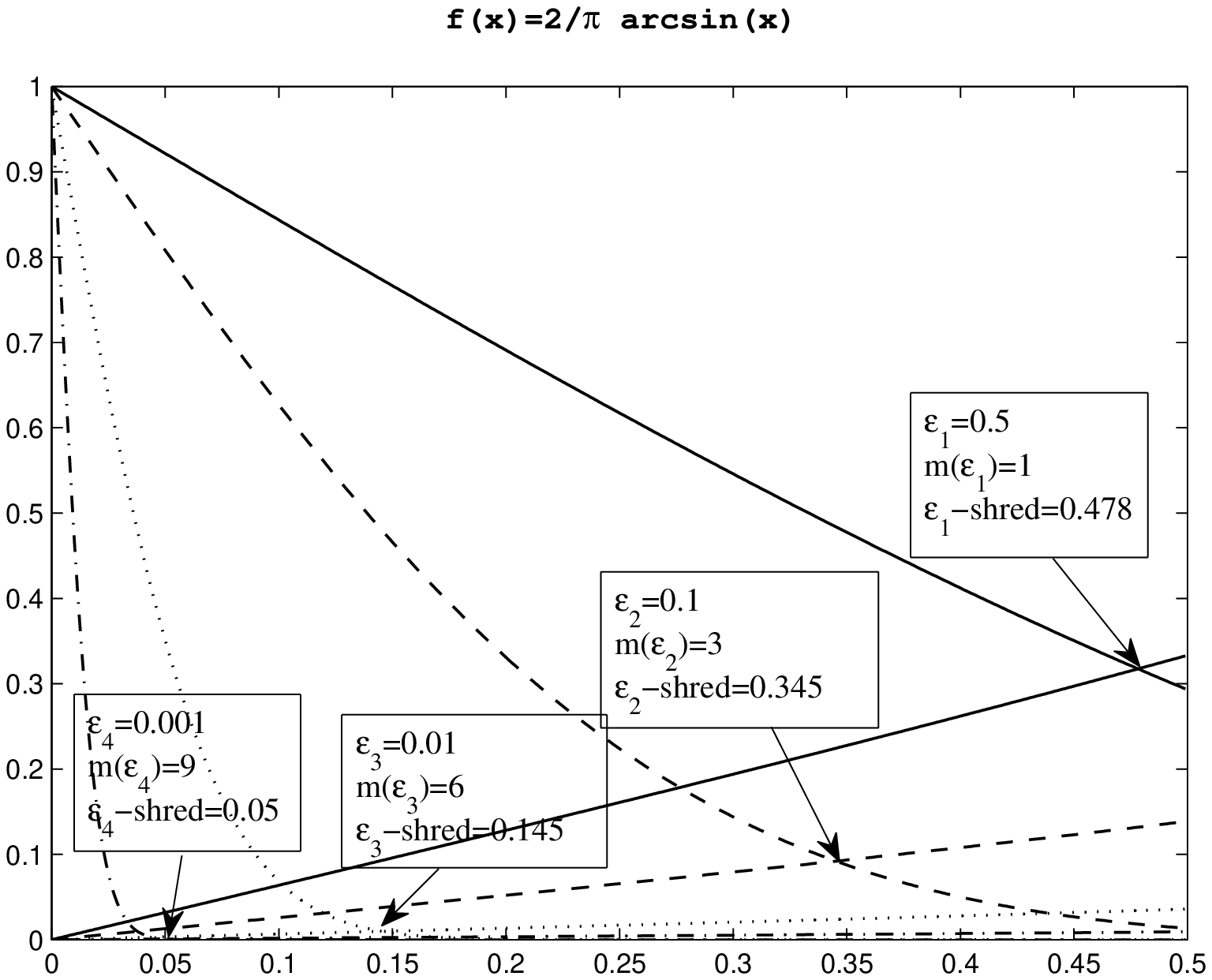}\\
 \caption{The $\tau_m(\varepsilon)$-functions and $\varepsilon$-basins-shreds for $f(x)=x^2$ (on the left, the $y$-axis cut to [0,0.35] for better clarity of the picture) and $f(x)=(2/\pi)\arcsin x$ (on the right, the $x$-axis cut to [0,0.5])}\label{wykresy}
\end{figure}

The interesting thing is the localization of $\varepsilon$-basins-shred in the given interval $(z^k,z^{k+1})$. We carried out the numerical simulations for the two functions, $f(x)=x^2$ and $f(x)=\frac{2}{\pi}\arcsin x$, $x\in[0,1]$, and $\varepsilon_1=0.5$, $\varepsilon_2=0.1$, $\varepsilon_3=0.01$ and $\varepsilon_4=0.001$. These functions properly extended onto $\mathbb{R}$, i.e. such that in every interval $[l,l+1]$ we have a copy of $f(x)$ on $[0,1]$ shifted $l$-units upward, induce orientation preserving circle homeomorphisms with a fixed point. In both cases the fixed point $x^{+}=0$ was attracting and the fixed point $x^{-}=1$ was repelling for $x\in (0,1)$. The results of this numerical experiment are presented in Fig.\ref{wykresy}, together with the graphs of $\tau^{+}_m$ and $\tau^{-}_m$, for $m=m(\varepsilon_1), m(\varepsilon_2), m(\varepsilon_3), m(\varepsilon_4)$. It seems that $\varepsilon$-basins shred $\widetilde{x}(\varepsilon)$ tends to $x^{-}$ as $\varepsilon\to 0$ for $f(x)=x^2$ and $\widetilde{x}(\varepsilon)\to x^{+}$ for $f(x)=\frac{2}{\pi}\arcsin x$. However, the $\varepsilon$-basins-shred for $\varepsilon\to 0$ probably could as well be an interior point of $(z^k,z^{k+1})$. Defining conditions on a function $f(x)$ which determine whether the $\varepsilon$-basins-shred tends to the repelling or to the attracting end-point of the interval or to some point in its interior, is an issue for further research.

\begin{theorem}\label{uniwersalne N_1}
Let $\varphi \colon S^1 \to S^1$ be a homeomorphism with a rotation number $\varrho(\varphi)=\frac{p}{q}$.

For $\varepsilon>0$ there exists
$N=N(\varepsilon)$ such that for every point $z\in S^1$ at least one of the following two conditions is satisfied:

  \begin{tabular}{p{0.8cm}p{13.5cm}}
  $\ref{uniwersalne N_1}.1$ &  there exists a periodic point $z_0^+\in \textrm{Per}(\varphi)$ such that $|\varphi^{nq+i}(z)-\varphi^i(z_0^+)|<\varepsilon$ for all $n\geq N$ and  $i=0,1,...,q-1$ \\
  or & \empty\\
  $\ref{uniwersalne N_1}.2$ & there exists a periodic point $z_0^-\in \textrm{Per}(\varphi)$ such that $|\varphi^{-nq-i}(z)-\varphi^{-i}(z_0^-)|<\varepsilon$ for all $n\geq N$ and  $i=0,1,...,q-1$,\\
 \end{tabular}\\
  i.e. after $Nq$ iterations forward or $Nq$ iterations backward we are always $\varepsilon$-close to one of the periodic orbits.
 \end{theorem}

\noindent{\bf Proof.} Assume firstly that $\varphi$ has $r$ different periodic points  $z^1, z^2,
\,\dots,\, z^r$. For each
$1\leq k \leq r$ we apply the properties $(iv)-(vi)$ with $z_0^-= z^k$
and $z_0^+=z^{k+1}$ (or $z_0^-= z^{k+1}$
and $z_0^+=z^{k}$ if $z^k$ is attracting and $z^{k+1}$ repelling for $z\in (z^k, z^{k+1})$). Set $N=\tilde{m}=\tilde{m}_k$. Then  $z\in B_k^+$ satisfies $\ref{uniwersalne N_1}.1$ and $z\in B_k^-$ satisfies \ref{uniwersalne N_1}.2. Every point $z\in (a_{\widetilde{m}}(\varepsilon), b_{\widetilde{m}}(\varepsilon))$ fulfills both \ref{uniwersalne N_1}.1 and \ref{uniwersalne N_1}.2.  Now, since $S^1= [z^1, z^2] \cup \, \dots\, [z^{k-1},z^k] \cup
[z^{r}, z^1]$, it is enough to take $N= \max_{1\leq
k\leq r+1} \, \tilde{m}_k$ (where $\tilde{m}_{r+1}$ corresponds to the interval
$[z^{r}, z^1]$) to get the statement.

Suppose now that $\#\textrm{Per}(\varphi)=\infty$. $\textrm{Per}(\varphi)$ is closed thus compact
subset of $S^1$. Fix  $\varepsilon>0$. The proof will be carried out in the following steps:

\begin{tabular}{p{0.3cm}p{13.8cm}}
  (1) & Let $z_0$ be a periodic point with the orbit $\mathcal{O}=\{z_0, z_1, \ldots, z_{q-1}\}$ for which there exists another periodic point $z_0^{\prime}$ with the orbit $\mathcal{O}^{\prime}=\{z_0^{\prime}, z_1^{\prime}, \ldots, z_{q-1}^{\prime}\}$ such that for every $i=0,1,\ldots, q-1$, $z_{i}$ and $z_i^{\prime}$ are consecutive periodic points, $z_i^{\prime}>z_i$ and for at least one $i_{*}\in\{0,1,\ldots,q-1\}$ we have $\vert z_i^{\prime}-z_i\vert\geq \varepsilon$. If it is not possible to find such a point $z_0$, then the distance between any two consecutive periodic points is smaller than $\varepsilon$ and the hypothesis of Theorem \ref{uniwersalne N_1} is satisfied in a trivial way. \\
  \end{tabular}

  \begin{tabular}{p{0.3cm}p{13.8cm}}
  (2) & Notice that the number of intervals $(z,z^{\prime})$ between consecutive periodic points $z$ and $z^{\prime}$ such that $\vert z^{\prime}-z\vert\geq \varepsilon$ is finite. Consequently, the number of pairs $\{\mathcal{O}, \mathcal{O}^{\prime}\}$ of periodic orbits $\mathcal{O}$ and $\mathcal{O}^{\prime}$ such as in (1) is finite. Denote as $D^{\varepsilon}$ the collection of all such ordered pairs $\{\mathcal{O}, \mathcal{O}^{\prime}\}$. \\
   \end{tabular}

  \begin{tabular}{p{0.3cm}p{13.8cm}}
  (3) & Let now $z$ be an arbitrary point on $S^1$.  If $z\in\textrm{Per}(\varphi)$, there is nothing to prove.
  If $z\not\in \textrm{Per}(\varphi)$ and $z$ does not lie in any interval $(z_{*},z_{*}^{\prime})$, where $z_*\in \mathcal{O}$ and $z_*^{\prime}\in\mathcal{O}^{\prime}$ with the pair of periodic  orbits $\{\mathcal{O}_*, \mathcal{O}_*^{\prime}\}\in D^{\varepsilon}$, then  every point of the full orbit $\{\varphi^i(z)\}_{i\in\mathbb{Z}}$ belongs to some interval between the two periodic points with length smaller than $\varepsilon$. As a result conditions \ref{uniwersalne N_1}.1 and \ref{uniwersalne N_1}.2 of Theorem \ref{uniwersalne N_1} are satisfied in a trivial way with arbitrary $N\in\mathbb{N}\cup\{0\}$. If $z\not\in\textrm{Per}(\varphi)$ but there exist periodic points $z_{*}$ and $z_*^{\prime}$ such that $z\in (z_*, z_*^{\prime})$ and $z_*\in \mathcal{O}_*$, $z_*^{\prime}\in\mathcal{O}_*^{\prime}$ with a pair of orbits $\{\mathcal{O}_*,\mathcal{O}_*^{\prime}\}\in D^{\varepsilon}$, then at least one of the conditions \ref{uniwersalne N_1}.1 or \ref{uniwersalne N_1}.2 holds for $z$ with $N=N_{\max}$, where $N_{\max}$ is the ``universal'' $N$ derived, as in the first part of the proof, for the finite collection of all the intervals $\{(z^i,z^{\prime \ i})\}$ between the two consecutive periodic points whose orbits $\mathcal{O}_i$ and $\mathcal{O}_i^{\prime}$ form pairs $\{\mathcal{O}_i,\mathcal{O}_i^{\prime}\} \in D^{\varepsilon}$.\\
  \end{tabular}

Consequently, it is enough to take $N=N_{\max}$ for an arbitrary $z\in S^1$.\hbx

Note that in general a number $N$ satisfying the statement of Theorem \ref{uniwersalne N_1} could be found by considering, instead of $\varepsilon$-basins-shreds, just the geometrical middles $\hat{z}_k$ of the intervals $\hat{I}_k$ between periodic points (with union $\bigcup_{k}\hat{I}_k$ giving the whole of $S^1$), computing the corresponding numbers $\hat{N}_k$ such that the iterates $\varphi^{-nq-i}(\hat{z}_k)$ and  $\varphi^{nq+i}(\hat{z}_k)$ are placed in the $\varepsilon$-neighbourhood of periodic orbits for $n>\hat{N}_k$ and $i=0,1,\dots,q-1$, and then taking $N=\max_k\hat{N}_k$. However, given $\varepsilon>0$ the notion of the $\varepsilon$-basins-shred says how to find the smallest $N\in\mathbb{N}$ with these properties.

\begin{proposition}\label{wlasnosc dla przemieszczenia} Let $\varrho(\varphi)=\frac{p}{q}$. Then for every $\varepsilon > 0$ there exists $N$ such that for every $z\in S^1$ the sequence
$\{\eta_n(z)\}_{n=-\infty}^{\infty}$ satisfies at least one of the following statements:

\begin{tabular}{p{1.5cm}p{13.5cm}}
  $\ref{wlasnosc dla przemieszczenia}.1$ &  $\forall_{n>N} \ \forall_{l\in \mathbb{N}} \ \vert\eta_{n+lq}(z)-\eta_n(z)\vert <\varepsilon$\\
  or & \empty\\
  $\ref{wlasnosc dla przemieszczenia}.2$ & $\forall_{n>N} \ \forall_{l\in\mathbb{N}} \ \vert\eta_{-(n+lq)}(z)-\eta_n(z)\vert <\varepsilon$\\
 \end{tabular}\\
\end{proposition}
\noindent{\textbf{Proof.}} The proposition is a direct consequence of Theorem \ref{uniwersalne N_1}. \hbx

\section{Homeomorphisms with irrational rotation number}

Let now $\varrho(\varphi)$ be irrational. If $\varphi$ is not transitive by $\Delta\subset S^1$ denote the unique minimal set of $\varphi$ (a Cantor type set), by $\widetilde{\Delta}\subset \mathbb{R}$  the total lift of  $\Delta$ to $\mathbb{R}$ and by $\widetilde{\Delta}_0\subset [0,1)$ the lift of $\Delta$ to $[0,1)$, i.e. $\widetilde{\Delta} = \{x+k; \ x\in\widetilde{\Delta}_0, \ k\in\mathbb{Z}\}$.

\begin{proposition}\label{density in the interval}
Let $\varphi: S^1\to\ S^1$  be a homeomorphism with the irrational rotation number $\varrho$ and $\Phi:\mathbb{R}\to\mathbb{R}$ its lift.

\begin{tabular}{p{0.8cm} p{13.3cm}}
  $\ref{density in the interval}.1$ & If $\varphi$ is transitive, then for every $z\in S^1$ the sequence $\eta_n(z)$, $n\in \mathbb{N}$,
is dense in the interval $[\min_{x\in [0,1]}\{\Phi(x)-x \ \rmod \ 1\}, \max_{x\in [0,1]}\{\Phi(x)-x \ \rmod \  1\}]$, which can be rewritten in the form\\
\end{tabular}
\begin{equation}\label{wzor na przedzial}
    \qquad \;\; \qquad \;\; [\min_{x\in[0,1]}\{\Gamma^{-1}(x+\varrho)-\Gamma^{-1}(x) \ \rmod \  1\},\max_{x\in[0,1]}\{\Gamma^{-1}(x+\varrho)-\Gamma^{-1}(x) \ \rmod \  1\}],
\end{equation}

\begin{tabular}{p{0.8cm}p{13.3cm}}
\empty & where  $\Gamma$ is a lift of a homeomorphism $\gamma$ conjugating $\varphi$ to the rotation $r_{\varrho}$,\newline $\;\;\;\;$ i.e. $\varphi=\gamma^{-1}\circ r_{\varrho}\circ \gamma$ \\

  $\ref{density in the interval}.2$ & If $\varphi$ is not transitive, then  \\ \end{tabular}\\
  \begin{tabular}{p{0.8cm}p{13.3cm}} \empty & \begin{tabular}{p{0.5cm}p{12.5cm}} i) & for  $z\in \Delta$ $\eta_n(z)$ is dense in a set \\ \end{tabular}\\
  \end{tabular}
\begin{equation}\label{wzor na zbior gestosci}
 \quad \quad \quad \quad \quad \quad \ \mathcal{D}:=\{\Phi(x)-x \ \rmod \ 1: \quad x\in \widetilde{\Delta}_0\}
\end{equation}
\begin{tabular}{p{0.8cm}p{13.3cm}} \empty & \begin{tabular}{p{0.5cm}p{12.5cm}} ii) & for $z\in S^1\setminus \Delta$ and $w\in \Delta$ there exist increasing sequences $\{n_k\}$ and $\{\widehat{n}_k\}$ such that for every $l\in\mathbb{Z}$ \\ \end{tabular}\\ \end{tabular}
$$\lim_{k\to\infty}\eta_l(\varphi^{n_k}(z))=\eta_l(w)\quad \textrm{and} \quad \lim_{k\to\infty}\eta_l(\varphi^{-\widehat{n}_k}(z))=\eta_l(w).$$
\end{proposition}
The statement \ref{density in the interval}.2 $ii)$ means that  $\mathcal{D}$ is $\omega$- and $\alpha$-limit set for displacements of points outside $\Delta$.

\noindent{\bf Proof.} For the proof of \ref{density in the interval}.1 fix $z_0\in S^1$ and consider the sequence $\eta_n(z_0)=\Phi^{n}(x_0)-\Phi^{n-1}(x_0) \ \rmod \  1$, $n\in \mathbb{N}$. Since $\Phi^{n-1}(x_0)=k +\alpha_{n-1}$ for some $k\in \mathbb{Z}$ and $\alpha_{n-1}  = \Phi^{n-1}(x_0) \ \rmod \  1$ being the fractional part of $\Phi^{n-1}(x_0)$,  $\Phi^{n}(x_0)=\Phi(k+\alpha_{n-1})=k+\Phi(\alpha_{n-1})$ and $\Phi^{n}(x_0) - \Phi^{n-1}(x_0)=\Phi(\alpha_{n-1})-\alpha_{n-1}.$ Hence we have the inclusion $$\{\Phi^{n}(x_0)-\Phi^{n-1}(x_0)\} \subset [\min_{x\in [0,1)} \Phi(x) -x, \max_{x\in [0,1)} \Phi(x)-x ].$$ Density of $\eta_n(z_0)$ in $[\min_{x\in[0,1]}\{\Phi(x)-x \ \rmod \  1\}, \max_{x\in[0,1]}\{\Phi(x)-x \ \rmod \  1\}]$ follows easily from  density of trajectories in $S^1$. We can write the interval of concentration also in the form of (\ref{wzor na przedzial}) since the displacement function  $\Psi(x)$ is periodic with period $T=1$, $\Gamma$ is a homeomorphism which maps the intervals of length $1$ onto the intervals of length $1$ and $\Gamma^{-1}(x+\varrho)-\Gamma^{-1}(x) =\Phi (\tilde{x})-\tilde{x}$, where $x=\Gamma(\tilde{x})$.

The proof of $\ref{density in the interval}.2 \ i)$ requires only minor modifications of the first part of the proof of $\ref{density in the interval}.1$. The statement $\ref{density in the interval}.2\ ii)$ follows from the fact that $\Delta$ is  $\omega$- and $\alpha-$ limit set for trajectories of points $z\in S^1\setminus \Delta$. \hbx

Since the conjugating homeomorphism is determined uniquely up to a rotation, i.e. up to an additive constant when we pass to a lift, the formula (\ref{wzor na przedzial}) involving the conjugating homeomorphism $\Gamma$ does not depend on the particular choice of $\gamma$ (or its lift $\Gamma$). For a rigid rotation $r_{\varrho}$ the ``interval of concentration'' is obviously the one-point set $\{\varrho\}$.

\subsection{Distribution of displacements}
A circle homeomorphism $\varphi$ with the irrational rotation number is always metrically isomorphic to the irrational rotation. The non-atomic unique $\varphi$-invariant Borel probability measure $\mu$ is given by $\mu(A)=\Lambda (\gamma(A)), \ A\subset S^1$, where $\Lambda$ denotes the normalized Lebesque'e measure on $S^1$ and $\gamma: S^1\to S^1$ is a continuous non-decreasing surjective map such that $\gamma\circ \varphi = r_{\varrho}\circ \gamma$. If $\varphi$ is transitive then $\gamma$ is a homeomorphism and thus $\varphi$ is conjugated to the rotation $r_{\varrho}$. Every $C^1$-diffeomorphism  with a derivative $\varphi^{\prime}$ of bounded variation is transitive on the account of Denjoy Thereom (\cite{denjoy}). Thus in particular every $C^2$-diffeomorphism is transitive.  If $\varphi$ is not transitive, then it is semi-conjugated to $r_{\varrho}$, the invariant measure $\mu$ is concentrated on the minimal set $\Delta$ (a Cantor-type set) and the semi-conjugacy $\gamma$ is $1\textrm{-to-}1$ on $\Delta$ with except some countable set $E\subset \Delta$, where $\mu(E)=0$. Moreover, $\gamma$ is constant at the intervals complementary to $\Delta$, equal to the value of $\gamma$ at the endpoints of a given complementary interval.

Firstly we will prove the following theorem, which is intuitively natural but not presented in accessible literature:
\begin{theorem}\label{twgs}
 The mapping $\varphi\mapsto\gamma$ assigning to a homeomorphism $\varphi$ with irrational rotation number $\varrho$ a map $\gamma:S^1\to S^1$ semiconjugating (or conjugating, if $\varphi$ is transitive) $\varphi$ with the rotation $r_{\varrho}$, is a continuous mapping from $C^{0}(S^1)$ into $C^{0}(S^1)$-topology.
\end{theorem}

 Theorem \ref{twgs} says that if $\varphi_1$ and $\varphi_2$ have irrational rotation numbers $\varrho_1$ and $\varrho_2$, respectively, and are close enough in $C^0(S^1)$, then the maps $\gamma_1$ and $\gamma_2$ (semi-)conjugating, correspondingly, $\varphi_1$ with $r_{\varrho_1}$ and $\varphi_2$ with $r_{\varrho_2}$, are arbitrarily close in $C^0(S^1)$ (possibly after the proper normalization). \\

\noindent{\bf Proof of Theorem \ref{twgs}.} Let $\varphi_1:S^1\to S^1$ be a homeomorphism with the rotation number $\varrho_1\in\mathbb{R}\setminus\mathbb{Q}$. Let $\gamma_1$ be a map that (semi-)conjugates $\varphi_1$ with the rotation $r_{\varrho_1}$. Denote by $\Delta_{\varphi_1}\subseteq S^1$ a minimal invariant set of $\varphi_1$. The orbit of an arbitrary point of $\Delta_{\varphi_1}$ is dense in $\Delta_{\varphi_1}$. Obviously, if  $\varphi_1$ is transitive then $\Delta_{\varphi_1}=S^1$.

Let $\varepsilon>0$ and $z_0\in\Delta_{\varphi_1}$ be arbitrary. There exists $\delta<\varepsilon/3$ such that $\vert \gamma_1(z_1)-\gamma_1(z_2)\vert<\varepsilon/4$ for any $z_1,z_2\in S^1$ where $\vert z_1-z_2\vert\leq \delta$. There also exists $N\in\mathbb{N}$ such that for every $z\in\Delta_{\varphi_1}$ one can find $i\in\{0,1,\dots, N\}$ satisfying $\vert\varphi_1^i(z_0)-z\vert<\delta/8$. Now, there exists a neighborhood $\Omega$ of $\varphi_1$ in $C^0(S^1)$ such that for every homeomorphism $\varphi_2\in\Omega$ with an irrational rotation number $\varrho_2$, transitive or not, we have $\vert\varphi_1^{i}(z)-\varphi_2^{i}(z)\vert<\delta/8$ and $\vert r_{\varrho_1}^i(z)-r_{\varrho_2}^i(z)\vert<\delta/8$ for every $z\in S^1$ and $i\in\{0,1,\dots,N\}$.

Assume now that $\varphi_2\in\Omega$ is (semi-)conjugated with $r_{\varrho_2}$ via $\gamma_2$. We can assume that $\gamma_1(z_0)=\gamma_2(z_0)=z_0$. Then $r_{\varrho_1}^{i}(z_0)=r_{\varrho_1}^{i}(\gamma_1(z_0))=\gamma_1(\varphi_1^{i}(z_0))$ and $r_{\varrho_2}^{i}(z_0)=r_{\varrho_2}^{i}(\gamma_2(z_0))=\gamma_2(\varphi_2^{i}(z_0))$. Consequently, $\vert \gamma_1(\varphi_1^{i}(z_0))-\gamma_2(\varphi_2^{i}(z_0))\vert=\vert r_{\varrho_1}^{i}(z_0)-r_{\varrho_2}^{i}(z_0)\vert <\delta/8$ and for $i\in\{0,1,\dots,N\}$ we can estimate:
\begin{eqnarray}\label{oszacowanie}
\vert \gamma_1(\varphi_2^i(z_0))-\gamma_2(\varphi_2^i(z_0))\vert &\leq& \vert \gamma_1(\varphi_2^i(z_0))-\gamma_1(\varphi_1^i(z_0))\vert +\vert \gamma_1(\varphi_1^i(z_0))-\gamma_2(\varphi_2^i(z_0))\vert\nonumber \\
  \empty &=& \vert \gamma_1(\varphi_2^i(z_0))-\gamma_1(\varphi_1^i(z_0))\vert+ \vert r^i_{\varrho_1}(z_0)-r^i_{\varrho_2}(z_0)\vert<\varepsilon/3
\end{eqnarray}
Since for every $z\in\Delta_{\varphi_1}$ there exists $i\in\{0,1,\dots,N\}$ such that $\vert\varphi_1^i(z_0)-z\vert<\delta/8$ and, simultaneously, $\vert\varphi_1(z_0)-\varphi_2(z_0)\vert<\delta/8$ for every $i\in\{0,1,\dots,N\}$, we obtain that for every $z\in\Delta_{\varphi_1}$ there exists $i\in\{0,1,\dots,N\}$ such that $\vert\varphi_2^i(z_0)-z\vert<\delta/4$.

Let us now take arbitrary $z\in\Delta_{\varphi_1}$. Then:

\begin{enumerate}[1.]
  \item if $z\in [\varphi_2^{i_1}(z_0),\varphi_2^{i_2}(z_0)]$ for some indexes $0\leq i_1,i_2\leq N$ ($i_1\neq i_2$), where $\vert \varphi_2^{i_1}(z_0)-\varphi_2^{i_2}(z_0)\vert\leq \delta$, then $\gamma_1(z)\in [\gamma_1(\varphi_2^{i_1}(z_0)),\gamma_1(\varphi_2^{i_2}(z_0))]$ with $\vert \gamma_1(\varphi_2^{i_1}(z_0))-\gamma_1(\varphi_2^{i_2}(z_0))\vert <\varepsilon/4$. As $\gamma_2(z)\in [\gamma_2(\varphi_2^{i_1}(z_0)),\gamma_2(\varphi_2^{i_2}(z_0))]$,  it follows from (\ref{oszacowanie}) that $\vert \gamma_1(z)-\gamma_2(z)\vert<\varepsilon$. Notice that at this point the proof is completed for the transitive case, since then $\Delta_{\varphi_1}=S^1$ and, as the orbits are dense in the whole $S^1$, we can always find indexes $0\leq i_1\neq i_2\leq N$ with such a property.
  \item if such an interval as above does not exist, then it means that for the point $z$ the distance between the two nearest neighbouring points, $\varphi_2^{i_1}(z_0)$ and $\varphi_2^{i_2}(z_0)$, of the $N$-orbit $\{\varphi_2^{i}(z_0)\}$, $i\in\{0,1,\dots,N\}$ exceeds $\delta$ (i.e. $z\in (\varphi_2^{i_1}(z_0),\varphi_2^{i_2}(z_0))$ with $\vert \varphi_2^{i_2}(z_0)-\varphi_2^{i_1}(z_0)\vert>\delta$). However, still one of these points, say $\varphi_2^{i_1}(z_0)$, is located within $\delta/4$ of $z$. Then $\vert\varphi_2^{i_2}(z_0)-z\vert>3\delta/4$. Consequently, the interval $(z+\delta/4,z+\delta/2)$ is complementary to $\Delta_{\varphi_1}$: If in this interval there were other points from $\Delta_{\varphi_1}$, this would again contradict the fact that every point of the minimal invariant set $\Delta_{\varphi_1}$ lies within the distance of $\delta/4$ of the $N$-orbit $\{\varphi_2^{i}(z_0)\}$, $i\in\{0,1,\dots,N\}$. Thus $\gamma_1(z+\delta/4)=\gamma_1(z+\delta/2)$, where by $z+\delta/4$ and $z+\delta/2$ we denote the points situated, respectively, $\delta/4$ and $\delta/2$ away from $z$ in the direction of $\varphi_2^{i_2}(z_0)$.
      \begin{enumerate}[2.1]
        \item Suppose that there are some points from the invariant set $\Delta_{\varphi_1}$ in the interval $[z+\delta/2,\varphi_2^{i_2}(z_0))$ and let $\widetilde{z}$ be the point which is closest to $z$  among all of these points. Then $\vert \varphi_{2}^{i_2}(z_0)-\widetilde{z}\vert<\delta/4$, because $\varphi_2^{i_2}(z_0)$ is the point of the $N$-orbit $\{\varphi_2^{i}(z_0)\}$, $i\in\{0,1,\dots,N\}$, located at the smallest distance to $\widetilde{z}\in\Delta_{\varphi_1}$. Then $\gamma_1(z)\in [\gamma_1(\varphi_2^{i_1}(z_0)),\gamma_1(\varphi_2^{i_2}(z_0))]$, where $\vert \gamma_1(\varphi_2^{i_2}(z_0))-\gamma_1(\varphi_2^{i_1}(z_0))\vert \leq \vert\gamma_1(\varphi_2^{i_2}(z_0))-\gamma_1(\widetilde{z})\vert+\vert\gamma_1(\widetilde{z})-\gamma_1(z+\delta/4)\vert+\vert\gamma_1(z+\delta/4)-\gamma_1(z)\vert+\vert\gamma_1(z)-\gamma_1(\varphi_2^{i_1}(z_0))\vert <  \varepsilon/4+0+\varepsilon/4+\varepsilon/4=3\varepsilon/4$.
Further, $\gamma_2(z)\in [\gamma_2(\varphi_2^{i_1}(z_0)),\gamma_2(\varphi_2^{i_2}(z_0))]$, where the endpoints of the interval $[\gamma_2(\varphi_2^{i_1}(z_0)),\gamma_2(\varphi_2^{i_2}(z_0))]$ are within $\varepsilon/3$ from the corresponding endpoints of the interval $[\gamma_1(\varphi_2^{i_1}(z_0)),\gamma_1(\varphi_2^{i_2}(z_0))]$. It follows that $\vert\gamma_2(z)-\gamma_1(z)\vert<2\varepsilon$.
        \item If $[z+\delta/2,\varphi_2^{i_2}(z_0))\in S^1\setminus\Delta_{\varphi_1}$, then $\gamma_1(z+\delta/2)=\gamma_1(\varphi_2^{i_2}(z_0))$ and $\vert \gamma_1(\varphi_2^{i_2}(z_0))-\gamma_1(\varphi_2^{i_1}(z_0))\vert \leq \vert\gamma_1(\varphi_2^{i_2}(z_0))-\gamma_1(z+\delta/2)\vert+\vert\gamma_1(z+\delta/2)-\gamma_1(z)\vert+
\vert\gamma_1(z)-\gamma_1(\varphi_2^{i_1}(z_0))\vert<\varepsilon/2$. Similarly as before we obtain that $\vert\gamma_2(z)-\gamma_1(z)\vert<2\varepsilon$.
      \end{enumerate}
\end{enumerate}
Hence for every $z\in\Delta_{\varphi_1}$ we have $\vert\gamma_2(z)-\gamma_1(z)\vert<2\varepsilon$

Let now $z\in S^1\setminus \Delta_{\varphi_1}$. Then
\begin{enumerate}[1.]
  \item if there exist points $z_1,z_2\in\Delta_{\varphi_1}$ such that $z\in (z_1,z_2)$ and $\vert z_2-z_1\vert\leq \delta$, then $\gamma_1(z)\in [\gamma_1(z_1),\gamma_1(z_2)]$, where $\vert \gamma_1(z_2)-\gamma_1(z_1)\vert<\varepsilon/4$, and  $\gamma_2(z)\in [\gamma_2(z_1),\gamma_2(z_2)]$, where $\vert \gamma_2(z_1)-\gamma_1(z_1)\vert<2\varepsilon$ and $\vert \gamma_2(z_2)-\gamma_1(z_2)\vert<2\varepsilon$ by what we have shown already. Consequently, $\vert\gamma_2(z)-\gamma_1(z)\vert<5\varepsilon$.
  \item otherwise for $z_1,z_2\in\Delta_{\varphi_1}$, where $z\in (z_1,z_2)$ and $z_1$ and $z_2$ are the two nearest neighbouring points of $\Delta_{\varphi_1}$ to $z$, we have $\vert z_2-z_1\vert>\delta$. Then the interval $(z_1,z_2)$ is complementary to $\Delta_{\varphi_1}$ and thus $\gamma_1(z_1)=\gamma_1(z_2)=\gamma_1(z)$. But again, $\vert\gamma_2(z_1)-\gamma_1(z_1)\vert<2\varepsilon$ and $\vert\gamma_2(z_2)-\gamma_1(z_2)\vert<2\varepsilon$ and thus $\vert\gamma_2(z)-\gamma_1(z)\vert<4\varepsilon$.
\end{enumerate}
We have shown that $\vert\gamma_1(z)-\gamma_2(z)\vert<5\varepsilon$ for arbitrary $z\in S^1$, which ends the proof.  \hbx

We will determine the distribution $\mu_{\Psi}$ of displacements for a homeomorphism $\varphi$ with respect to the measure $\mu$ (by distribution we mean a normalized measure). The distribution will be given on $[0,1)$ since equivalently one can consider a system $([0,1], T)$ with $T=\Phi \ \rmod \  1$. Then $\Gamma(x)=\mu([0,x])$, $x\in\mathbb{R}$ and $([0,1], G)$, $G=\Gamma \ \rmod \  1$ is a (semi-)conjugating system. However, if we choose  the lifts $\Phi\colon \mathbb{R}\to\mathbb{R}$ and $\Gamma\colon \mathbb{R}\to\mathbb{R}$ such that $\Phi(0)\in (0,1)$ and $\Gamma(\Phi(x))=\Gamma(x)+\varrho$ ($\varrho\in (0,1)$) we can skip ``$\rmod \  1$'' and identify $\Phi\sim T$ and $\Gamma\sim G$ in the remaining part of this subsection. In case $\varphi$ is not transitive, by $\widehat{\Delta}$ we denote the lift of the set $\Delta\setminus E\subset S^1$ to $[0,1)$ and by  $\widehat{\Gamma}$ the lift $\Gamma$ cut to $\widehat{\Delta}$, i.e. $\widehat{\Gamma}=\Gamma\upharpoonright_{\widehat{\Delta}}$.

\begin{proposition}\label{twosupp} Let $\mu$ be the unique normalized invariant ergodic measure for a homeomorphism $\varphi:S^1 \to S^1$.

If $\varphi$ is transitive then the distribution $\mu_{\Psi}$ is the transported measure:
\begin{equation}\label{zwiazek nipsi lambda}
\mu_{\Psi}(A)\colon = \mu(\{x\in [0,1): \ \Psi(x)\in A\})=\Omega_*\Lambda(A)\colon=\Lambda(\Omega^{-1}(A)),
\end{equation}
where $\Omega(x):=\Gamma^{-1}(x+\varrho)-\Gamma^{-1}(x)$. The support of the displacement measure equals
\begin{equation}\label{nosnik}
\textrm{supp}(\mu_{\Psi})=\Psi([0,1])=\Omega([0,1]).\nonumber
\end{equation}

If $\varphi$ is not transitive then
\begin{equation}\label{drugizwiazek}
\mu_{\Psi}(A)=\Lambda(\widehat{\Omega}^{-1}(A)),
  \end{equation}
where $\widehat{\Omega}(x):=\widehat{\Gamma}^{-1}(x+\varrho)-\widehat{\Gamma}^{-1}(x)$ and
\begin{equation}
\textrm{supp}(\mu_{\Psi})=\Psi(\widetilde{\Delta}_0).\nonumber
\end{equation}
\end{proposition}
\noindent{\bf Proof.} Assume firstly that $\varphi$ is transitive. Let $A\subset [0,1)$. Then
\begin{eqnarray}
 \mu_{\Psi}(A)=\mu(\{x\in [0,1]: \ \Psi(x)\in A\})
 =\mu(\{x\in [0,1]: \ \Gamma^{-1}(\Gamma(x)+\varrho)-\Gamma^{-1}(\Gamma(x))\in A\}) \nonumber \\
 = \Lambda(\{x\in [0,1]: \Gamma^{-1}(x+\varrho)-\Gamma^{-1}(x)\in A\})=\Lambda(\Omega^{-1}(A)). \nonumber
\end{eqnarray}
 If $\varphi$ is not transitive then we perform similar calculations as above to obtain (\ref{drugizwiazek}). Supports of distributions are derived easily, too. \hbx

 Let now $\varphi_n$ be a sequence of homeomorphisms with irrational rotation numbers converging to $\varphi$ in the topology of $C^0(S^1)$ and $\mu^{(n)}$ a sequence of normalized $\varphi_n$-invariant ergodic measures. Let $\mu^{(n)}_{\Psi_n}$ be the corresponding displacement distributions. Before we will show that $\mu^{(n)}_{\Psi_n}$ converges weakly to $\mu_{\Psi}$, we recall some definitions and facts:

\begin{definition}{\bf (see e.g. \cite{bartoszynski})} Let $X$ be a complete separable metric space and $\mathcal{M}(X)$ the space of all finite measures defined on the Borel $\sigma$-field $\mathcal{B}(X)$ of subsets of $X$.

A sequence $\mu_n$ of elements of $\mathcal{M}$ is called \emph{weakly convergent} to $\mu\in\mathcal{M}(X)$ if for every bounded and continuous function $f$ on $X$
\begin{equation}
\lim_{n\to\infty} \int\limits_{X}f(x)\, d\mu_n(x)=\int\limits_{X}f(x)\, d\mu(x).
\end{equation}
We denote the weak convergence as $\mu_{n} \Longrightarrow \mu$.
\end{definition}
\begin{definition} A Borel set $A$  is said to be a \emph{continuity
set} for $\mu$ if A has $\mu$-null boundary, i.e. $\mu(\partial A)=0$.
\end{definition}
One can show (cf. \cite{rao}) that  $\mu_n\Longrightarrow \mu$ if and only if for each continuity set $A$ of $\mu$
$\lim_{n\to\infty}\mu_n(A)=\mu(A).$

\begingroup
\setcounter{tmp}{\value{theorem}}
\setcounter{theorem}{0} 
\renewcommand\thetheorem{\Alph{theorem}}
\begin{theorem}\label{twa}{\bf (cf. \cite{rao})} Suppose that $\mu_{n}\Longrightarrow \mu$. If $\xi^{(n)}$ is a sequence of continuous maps of $X$ into a metric space $Y$, converging uniformly on compacta to $\xi$, then $\xi_*^{(n)}\mu_n\Longrightarrow \xi_*\mu$.
\end{theorem}
\endgroup

\setcounter{theorem}{\thetmp} 

\begin{proposition}\label{pierwszazbieznosc} Suppose that $\varphi_n$ is a sequence of homeomorphisms with irrational rotation numbers $\varrho_n$ which converges in the metric of $C^0(S^1)$ to the homeomorphism $\varphi$ with irrational rotation number $\varrho$. Let $\mu^{(n)}_{\Psi_n}$ and $\mu_{\Psi}$ be the corresponding displacement distributions with respect to the invariant measures  $\mu^{(n)}$, $\mu$. Then
$\mu^{(n)}_{\Psi_n}\Longrightarrow \mu_{\Psi}$.
\end{proposition}
\noindent{\bf Proof.}  Since $\mu^{(n)}=\gamma^{(n)}_*\Lambda$ and $\mu=\gamma_*\Lambda$, where $\gamma^{(n)}$ are the maps (semi-)conjugating $\varphi_n$ with the corresponding rotations, Theorem \ref{twgs} and Theorem \ref{twa} yield that $\mu^{(n)}\Longrightarrow \mu$. But then also $\mu^{(n)}_{\Psi_n}\Longrightarrow \mu_{\Psi}$ on the account of Theorem \ref{twa}, because the displacement distributions are the invariant measures transported by the displacement functions $\Psi_n\to \Psi$ in $C^0(\mathbb{R})$. \hbx

For the ``good'' properties of the invariant measure and the distribution of displacements we need to assure that the conjugacy $\gamma$ is a $C^1$-diffeomorphism. The sufficient conditions for this are given, for example,  by Theorem \ref{twb} (the celebrated theorem of Herman, cf. eg. \cite{herman} or \cite{khanin sinai}), Yoccoz Theorem (\cite{yoccoz}) or Theorem \ref{twc} (\cite{khanin tep}):
\begingroup
\setcounter{tmp}{\value{theorem}}
\setcounter{theorem}{1} 
\renewcommand\thetheorem{\Alph{theorem}}
\begin{theorem}\label{twb}
Let $\varphi:S^1\to S^1$ be a diffeomorphism with irrational rotation number $\varrho$ which is conjugated to the rotation $r_{\varrho}$ by a homeomorphism $\gamma:S^1\to S^1$. Suppose that the following two conditions are satisfied
\begin{description}
  \item[B.1] $\varphi\in C^{2+\nu}$, $\nu>0$, $\varphi^{\prime}\geq \textrm{const}>0$;
  \item[B.2] if $\varrho=[k_1,k_2,\ldots, k_n, \ldots]$ is the expansion of $\varrho$ into continued fraction, then $k_n\leq n^{\kappa}$, $\kappa>0$.
\end{description}
Then $\gamma$ is a $C^1$-diffeomorphism.
\end{theorem}

\begin{theorem}\label{twc} Let $\varphi$ be a $C^{2+\alpha}$-smooth  circle diffeomorphism with rotation number $\varrho$ in Diophantine class $D_{\delta}$, where $0<\delta<\alpha\leq 1$ (for the definition of Diophantine class $D_{\delta}$ see \cite{khanin tep}).

Then the conjugacy $\gamma$ is  $C^{1+\alpha-\delta}$-smooth.
\end{theorem}
\endgroup

\setcounter{theorem}{\thetmp} 

The uniform convergence of $\mu^{(n)}_{\Psi_n}$ to $\mu_{\Psi}$ is rarely to occur, even just on some particular subclass of Borel sets. With the use of the theorem below we will see what kind of assumptions can guarantee the uniform convergence of displacement distributions:
\begingroup
\setcounter{tmp}{\value{theorem}}
\setcounter{theorem}{3} 
\renewcommand\thetheorem{\Alph{theorem}}
\begin{theorem}\label{twd}
{\bf (cf. \cite{rao})} Suppose $\mu$ is a measure on $\mathbb{R}^k$ such that every convex subset of $\mathbb{R}^k$ is a continuity set for $\mu$. Then $\mu_n\Longrightarrow\mu$ if and only if
\begin{equation}\label{sup}
\sup[\vert \mu_n(C)-\mu(C)\vert, \ C\in\mathcal{C}]\to 0,
\end{equation}
where $\mathcal{C}$ denotes the class of all measurable convex sets.
\end{theorem}

\endgroup

\setcounter{theorem}{\thetmp} 

 In particular (\ref{sup}) is true when $\mu$ is absolutely continuous with respect to the Lebesque measure. We will obtain a similar result as (\ref{sup}) for the displacement distribution $\mu_{\Psi}$. However, note that even if the invariant measure $\mu$ is absolutely continuous, it might not be so for the displacement distribution, as happens for example for the rotation $r_{\varrho}$, where $\mu_{\Psi}$ is the singular measure equal to the Dirac delta $\delta_{\varrho}$.
\begin{proposition}\label{jednostajnaslaba}  Let  $\varphi_n$ be a sequence of homeomorphisms with irrational rotation numbers converging in $C^0(S^1)$ to a diffeomorphism $\varphi$ with an irrational rotation number $\varrho$. Suppose  that $\gamma$ conjugating $\varphi$ with $r_{\varrho}$ is a $C^1$-diffeomorphism and that the set of critical points of the displacement function $\Psi(x) = \Phi(x) - x $, i.e. the set $C:=\{x\in \mathbb{R}: \ \Phi^{\prime}(x)=1\}$, is of Lebesque measure $0$. Then
\begin{equation}\label{jed} \sup[\vert \mu^{(n)}_{\Psi_n}(I)-\mu_{\Psi}(I)\vert, \ I\in \mathcal{J}]\to 0,
\end{equation}
where $\mathcal{J}$ denotes the collection of all the intervals $I\subset [0,1]$ (open, closed or half-closed).
\end{proposition}
\begin{lemma}\label{absolutely continuous} Under the assumptions of Proposition \ref{jednostajnaslaba} the distribution $\mu_{\Psi}$ is absolutely continuous with respect to the Lebesque measure. \end{lemma}
\noindent{\bf Proof.}
Notice that the fact that $\gamma$ is a $C^1$-diffeomorphisms is equivalent to the absolute continuity of the measure $\mu$ with the continuous density function. Suppose then that $\Lambda(B)=0$ and let $A=\Psi^{-1}(B)$, $A_1=A\cap C$ and $A_2=A\cap (\mathbb{R}\setminus C)$. By the Inverse Function Theorem, for every point $x\in A_2\subset \mathbb{R}\setminus C$ there exists an open neighbourhood $\mathcal{B}(x,r_x)$ such that $\Psi\upharpoonright_{\mathcal{B}(x,r_x)}$ is a $C^1$-diffeomorphism. Let $A_x=A_2\cap \mathcal{B}(x,r_x)$. Since every diffeomorphism preserves sets of Lebesque measure $0$, $\Lambda(A_x)=\Lambda(B_x)=0$, where $B_x\colon =\Psi(A_x)\subset B$. Finally, from the cover $\{A_x\}_{x\in A_2}$ we can choose a countable subcover, because $A_2$ is separable. Therefore $\Lambda(A_2)=0$.  But since $\mu$ has a density, $\Lambda(A)=\Lambda(A_2)=0$ already implies $\mu_{\Psi}(B)=\mu(A)=0$.\hbx

We will give an effective description of the density function of
$\mu_\Psi$ with respect to  the Lebesque measure $\Lambda$.

Denote by $V_\Psi(C)$ the set of critical values of $\Psi$. The following statement can be justified with the use of the Inverse Function Theorem:
\begin{lemma}\label{finite counterimage}
Let $\Psi = \Phi - {\rm I}: \mathbb{R} \to \mathbb{R}$, where $\Phi$ is a
lift of a $C^1$-diffeomorphisms $\varphi:S^1\to S^1$.

Then for every $y\in \mathbb{R} \setminus V_\Psi(C)$ the set
$\Psi^{-1}(y)\cap [0,1]$ is finite.
\end{lemma}

\begin{theorem}\label{effective density} Let $\Psi= \Phi - {\rm I}$  be a displacement function of  a lift $\Phi$ of  a transitive $C^1$-diffeomorphism $\varphi$ of the circle. Assume next that  $\Phi$  is conjugated with the translation by $\varrho(\varphi)$
by  $C^1$-diffeomorphism $\Gamma$ and that the set $C_\Psi\subset
[0,1]$  of critical points of $\Psi$ is of  Lebesque measure
equal to $0$.

Then  the
 density  function  $\Delta(y)$ of $\mu_\Psi$ with respect to the Lebesgue measure $\Lambda$  is  equal to
 $$ \Delta(y)=\left\{
  \begin{array}{ll}
    0 & \hbox{if $y\notin {\rm supp}(\mu_\Psi)$;} \\
   \sum\limits_{x\in \Psi^{-1}(y)} \Gamma^{\prime}(x)\, \frac{1}{\vert\Phi^{\prime}(x)-1\vert} & \hbox{if $y\in {\rm supp}(\mu_\Psi)$.}
  \end{array}
\right.
  $$
where the latter is well-defined almost everywhere in ${\rm
supp}(\mu_\Psi)$, i.e. in  ${\rm supp}(\mu_\Psi) \setminus
V_\Psi(C)$.
\end{theorem}
\noindent{\bf Proof.} Note that in this case ${\rm supp}(\mu_\Psi)$ is a non-degenerate interval in $[0,1]$ as follows from Proposition \ref{twosupp}. Let then ${\rm supp}(\mu_\Psi) = [a,b] \subset [0,1]$. We show that the function  $F(y)=\mu_\Psi([a, y])$ is strictly monotonic.
 Indeed, if $y_1 < y_2$ then $F(y_2)- F(y_1)= \mu_\Psi((y_1, y_2])$. It is enough to show that $\mu_\Psi((y_1, y_2])>0$. Let $y\in (y_1, y_2] $ be a regular value of $\Psi$ and $\Psi^{-1}(y)= \{x_1, \, \dots\,,x_n\}$. There exist
 open  balls $B_i=B(x_i,r_i)$  such that $\Psi\upharpoonright_{B(x_i,r_i)}$ is a
 diffeomorphism onto an open neighborhood $U_i = \Psi(B(x_i,r_i)) \subset (y_1,y_2)$  of
 $y$. Let $U= \bigcap_{i=1}^n U_i$ and $ V_i =
 \Psi^{-1}(U)\cap B_i$.  We have $\mu_\Psi((y_1,y_2]) \geq \mu_\Psi(U) = \mu(\Psi^{-1}(U))=\mu(\bigcup_{i=1}^n V_i)=\sum_{i=1}^n \mu(V_i)$. On the other hand $\mu$ has a continuous density $\Gamma^\prime$
 with respect to the Lebesgue measure $\Lambda$.  This implies that
 $\mu(V_i)>0$ for every $1\leq i \leq n$,  because $\Gamma^\prime(x)$ cannot be identically
 equal to $0$ on an open set $V_i$ as the derivative of  a
 diffeomorphism.

 Since $F(y)$ is strictly monotonic on $(a,b]$ it has a measurable  derivative $f(y)$  almost everywhere in $[a,b]$ with respect to the Lebesgue measure. Moreover, $f$ is integrable and consequently
  we have $F(y)= \int_a^y f(x)d x$. Furthermore, any two such derivatives are equal up to a set of the Lebesgue  measure zero. To  complete the proof it is enough to derive the derivative of $F(y)$ at a point $y\in {\rm supp}(\mu_\Psi) \setminus
V_\Psi(C)$.

To derive this derivative, first observe that $F(y+h)- F(y)=
\mu_\Psi((y, y+h]) =\mu_\Psi([y,y+h])$ (for $h>0$, but for $h<0$ everything what follows is analogous), because $\mu_\Psi$ is
absolutely continuous with respect to the Lebesgue measure (Lemma
\ref{absolutely continuous}). For $h$ sufficiently small  take $\Psi^{-1}([y, y+h])
\subset [0,1]$.  The set $C_\Psi\subset [0,1]$ is closed, thus
compact, subset of Lebesgue measure $0$ by our assumption. Therefore also $\mu(C_\Psi)=0$. Put $R=
[0,1]\setminus C_\Psi$.

Now we consider another measure $\tilde{\mu}_\Psi$ defined by
the formula
$$ \tilde{\mu}_\Psi(A) = \mu(\Psi^{-1}(A)\cap R)$$
By our assumption on $C_\Psi$, we have
$\mu_\Psi(A)=\tilde{\mu}_{\Psi}(A)$ for every measurable set $A \subset
[0,1]$. Consequently, it is sufficient to derive the derivative of
$\tilde{F}(y)= \tilde{\mu}([a,y])$ at a regular value $y \in
\textrm{supp}\mu_{\Psi}$.

Let then $y$ be a regular value. Observe that for $h>0$ small enough the entire interval $[y,y+h]$ is contained in
$\Psi(R)$ and then consequently  it consists of regular values of
$\Psi$. Suppose that $\Psi^{-1}(y)=\{x_1, \, \dots\,, x_n\}\subset [0,1]$ (by Lemma \ref{finite counterimage} the set $\Phi^{-1}(y)\cap [0,1]$ is finite). Consider the set $B = \Psi^{-1}([y, y+h])\cap [0,1]$. It is a closed, thus compact set in [0,1] and  it contains at least $n$
connected components $\mathcal{C}_1, \,\dots\, , \mathcal{C}_n$
containing $x_i$, $1\leq i\leq n$ respectively. Each
$\mathcal{C}_i$ is an interval contained in $R$. Possibly the number of connected components in $B$ could be even infinite but one can show that for sufficiently small $h>0$ $B$ consists of exactly $n$ connected components $\mathcal{C}_i=\mathcal{C}_i(h)$, i.e. $B=\bigcup_{i=1}^{n}\mathcal{C}_i$.

Now our task is to estimate the measure of $$\tilde{\mu}_\Psi([y,y+h])= \mu_\Psi([y,y+h])=\sum_{i=1}^n \mu(\mathcal{C}_i)\,.$$ Since $\mu$ is absolutely
continuous with respect to  $\Lambda$ with the density
$\Gamma^\prime$, and $\Psi^{-1} : [y, y+h] \to \mathcal{C}_i$ is a
diffeomorphism, we have
$$ \mu(\mathcal{C}_i) = \int_{\mathcal{C}_i} \, \Gamma^\prime(x)
d \Lambda(x) = \int_y^{y+h} \Gamma^\prime(\Psi^{-1}_{(i)}(t))\vert
(\Psi^{-1}_{(i)})^\prime(t) \vert d t,$$
where $\Psi^{-1}_{(i)}$ denotes the diffeomorphism $\Psi^{-1}: [y,y+h]\to \mathcal{C}_i$.

Finally, by the Mean Value Theorem, we can replace the last
integral by
$$\Gamma^\prime(\Psi^{-1}_{(i)}(\theta_i))\vert
(\Psi^{-1}_{(i)})^\prime(\theta_i) \vert h $$
 where $ \theta_i \in [y, y+h]$. By continuity of $\Gamma^\prime $ and
 $(\Psi^{-1}_{(i)})^\prime$ this expression divided by $h$  tends to
 $$\Gamma^\prime(\Psi^{-1}_{(i)}(y)) \, \vert (\Psi^{-1}_{(i)})^\prime(y)\vert = \Gamma^{\prime}(x_i)\frac{1}{\vert \Phi^{\prime}(x_i)-1\vert}$$
 which gives the required formula for the density $\Delta$  of
 $\mu_\Psi$. \hbx

\noindent{\bf Proof of Proposition \ref{jednostajnaslaba}.} As $\mu_{\Psi}\ll \Lambda$, $\mu_{\Psi}(\partial A)=0$, where $A=I$ is an arbitrary interval. On the account of Proposition \ref{pierwszazbieznosc}  and Theorem \ref{twd} the convergence of measures $\mu^{(n)}_{\Psi_n}$ to $\mu_{\Psi}$ is uniform on the collection $\mathcal{J}$ of all intervals in $[0,1)$. \hbx

Note that the assumption of Proposition \ref{jednostajnaslaba} that the conjugacy of $\varphi$ with the rotation is of class $C^1$, i.e. that the unique invariant Borel probability measure is absolutely continuous with respect to the Lebesque measure, excludes critical circle homeomorphisms for which the invariant probability is always singular (\cite{graczykswiatek}). Moreover, circle diffeomorphisms having break points (discontinuities in the derivative, see e.g. \cite{dzahilov}) and even analytic circle diffeomorphisms (see \cite{arnold} for a particular example) might have singular invariant measures.

It is also worth emphasizing that in order to prove Proposition  \ref{jednostajnaslaba}, Lemma \ref{absolutely continuous} and Theorem \ref{effective density} we needed not only continuous differentiability of $\gamma$ but also the assumption on the measure of the set of critical points of $\Psi$.

Observe that the displacements distribution is ``invariant along the orbits'' i.e. when one defines $\Psi_{(1)}(x)=\Phi(x)-x$, $\Psi_{(2)}(x)=\Phi^2(x)-\Phi(x)$, ..., $\Psi_{(n)}(x)=\Psi(\Phi^{n-1}(x))=\Phi^n(x)-\Phi^{n-1}(x)$, $n\in\mathbb{N}$, and for the negative semi-orbit $\Psi_{(0)}(x)=x-\Phi^{-1}(x)$, $\Psi_{(-1)}(x)=\Phi^{-1}(x)-\Phi^{-2}(x)$, ..., $\Psi_{(-n)}(x)=\Psi(\Phi^{-(n+1)}(x))=\Phi^{-n}(x)-\Phi^{-(n+1)}(x)$, then using that fact the $\mu$ is $\varphi$ invariant we obtain
\begin{remark} For every $m,n\in \mathbb{Z}$ and $A\subset [0,1]$
$$\mu_{\Psi_{(n)}}(A)=\mu_{\Psi_{(m)}}(A)$$
\end{remark}

Let us make the easy observation  that the measure $\mu_{\Psi}(A)$ can be approximated by measuring the average frequency of points $\Phi^i(x)$ with values $\Psi(\Phi^i(x))$ in $A$ along a trajectory $\{\Phi^i(x)\}$, $i\in\mathbb{N}$:
\begin{proposition}\label{birkof} Let $\varphi: S^1\to S^1$ be a homeomorphism with the irrational rotation number $\varrho$ and $A\subset [0,1)$. Then
$$\frac{\sharp \{0\leq i \leq n-1: \ \Psi(\Phi^i(x))\in A\}}{n}\to\mu_{\Psi}(A)
$$
uniformly with respect to $x\in [0,1]$ as $n\to\infty$.
\end{proposition}
\noindent{\bf Proof.} The statement follows from the Birkhoff Ergodic Theorem applied for  an observable $f=\chi_{A}\circ \Psi$ and the uniformity  of the convergence over $x\in[0,1)$ follows from unique ergodicity of $\varphi$ (cf. Theorem \ref{wlasnosc uniquely ergodic}). \hbx

Naturally, also
\begin{remark}\label{srednia}
The rotation number $\varrho$ is the average of the distribution $\mu_{\Psi}$:
$$\int\limits_{[0,1)}\Psi d\mu = \varrho$$
\end{remark}

We will see that measuring displacements along the orbits  of a homeomorphism $\widetilde{\varphi}$, which is close enough to $\varphi$, also gives the approximation of the distribution $\mu_{\Psi}$ of $\varphi$. For $n\in\mathbb{N}$ and $z=\exp(2\pi\rmi x)\in S^1$ define a sample displacements distribution:
\[
\omega_{n,x}=\frac{1}{n}\sum_{i=0}^{n-1}\delta_{\Psi(\Phi^i(x))},
\]
where $\delta_{y}(A)=1$ if $y\in A$ and $\delta_{y}(A)=0$ otherwise, $A\subset [0,1]$. The Fortet-Mourier metric (see, for example, \cite{fortet-mourier})
is a method to measure the distance between the two probability measures on a given space:
\begin{definition}
Let $\mu$ and $\nu$ be the two Borel probability measures on a measurable space $(\Omega, \mathcal{F})$, where $\Omega$ is a compact metric space. Then the distance between the measures $\mu$ and $\nu$ is defined as
\[
d_{F}(\mu,\nu):=\sup\{\vert \int\limits_{\Omega} f\;d\mu - \int\limits_{\Omega} f\;d\nu\vert: \ f \ \textrm{is $1$-Lipshitz}\}\]
\end{definition}

The Fortet-Mourier metric is sometimes called also the  Wasserstein metric. It metrizes weak convergence of probability measures on spaces of bounded diameter (cf. eg. \cite{probmetric}). For the proof of the following theorem see \cite{fortet-mourier}:
\begingroup
\setcounter{tmp}{\value{theorem}}
\setcounter{theorem}{4} 
\renewcommand\thetheorem{\Alph{theorem}}
\begin{theorem}\label{tw o zbieznosci}
Assume that $P$ is a Borel probability on $(\Omega, \mathcal{F})$, which is a compact metric space, and let $(\omega_i)_{N\in\mathbb{N}}$ be a sample. Let $P_N:=\frac{1}{N}\sum_{i=1}^{N}\delta_{\omega_i}$ be the associated empirical law of $P$. Suppose that the support of $P$ is closed. Then
\[
\lim_{N\to\infty}d_F(P,P_N)=0
\]
\end{theorem}

\endgroup

\setcounter{theorem}{\thetmp} 
Let $\widetilde{\varphi}:S^1\to S^1$ be another homeomorphism, with rational or irrational rotation number, and let $\widetilde{\Phi}$ and $\widetilde{\Psi}=\widetilde{\Phi}-\textrm{Id}$ be a lift  and the displacement function of $\widetilde{\varphi}$, respectively. Define a sample displacements distribution associated with $\widetilde{\varphi}$: $\widetilde{\omega}_{n,x}:=\frac{1}{n}\sum_{i=0}^{n-1}\delta_{\widetilde{\Psi}(\widetilde{\Phi}^i(x))}$.

\begin{theorem}\label{przyblizanie wymiernymi}
Let $\varphi$ be a homeomorphism with irrational rotation number and the displacement distribution $\mu_{\Psi}$ with respect to the invariant measure $\mu$. For every $\varepsilon>0$ there exists a neighborhood $\mathcal{U}\subset C^{0}(S^1)$ of $\varphi$ such that for every homeomorphism $\widetilde{\varphi}\in \mathcal{U}$ and every $x_0\in[0,1]$ we have
\[
d_F(\lim_{n\to\infty}\frac{1}{n}\sum_{i=0}^{n-1}\delta_{\widetilde{\Psi}(\widetilde{\Phi}^i(x_0))},\mu_{\Psi})<\varepsilon
\]
\end{theorem}
In the proof we will make use of the fact that every circle homeomorphism with irrational rotation number (either transitive or not) is uniquely ergodic and of the following classical result (see, for instance, \cite{walters}):
\begingroup
\setcounter{tmp}{\value{theorem}}
\setcounter{theorem}{5} 
\renewcommand\thetheorem{\Alph{theorem}}
\begin{theorem}\label{wlasnosc uniquely ergodic}
If $T:X\to X$ ($X$-metrizable compact space) is uniquely ergodic then for every continuous function $f$ the time averages $\frac{1}{n}\sum_{i=0}^{n-1}f(T^i(x))$ converge uniformly.
\end{theorem}

\endgroup

\setcounter{theorem}{\thetmp} 
\noindent{\textbf{Proof of Theorem \ref{przyblizanie wymiernymi}.}} Fix $\varepsilon>0$. Notice that $supp(\mu_{\Psi})$ is closed, regardless of whether $\varphi$ is transitive or not. Thus on the account of Theorem \ref{tw o zbieznosci} for every $x\in[0,1]$ we have
\[\lim_{n\to\infty}d_F(\omega_{n,x},\mu_{\Psi})=0\]
However, since for every continuous function $f$ there is a uniform convergence
\[\int\limits_{[0,1]}fd\omega_{n,x}=\frac{1}{n}\sum_{i=0}^{n-1}f(\Psi(\Phi^i(x)))\to \int\limits_{[0,1]}f\circ \Psi d\mu=\int\limits_{[0,1]}f d\mu_{\Psi},\]
there exists $N\in\mathbb{N}$ such that  for all $n\geq N$ and $x\in[0,1]$ $d_{F}(\omega_{n,x},\mu_{\Psi})<\varepsilon/2$.
Further, there exists a neighbourhood $\mathcal{U}$ of $\varphi$ in $C^{0}(S^1)$ such that if $\widetilde{\varphi}\in \mathcal{U}$ is another homeomorphism with a lift $\widetilde{\Phi}$, then $\vert\widetilde{\Phi}^i(x)-\Phi^{i}(x)\vert<\varepsilon/4$ for every $0\leq i\leq N$ and $x\in[0,1]$. If $\widetilde{\omega}_{N,x}=\frac{1}{N}\sum_{i=0}^{N-1}\delta_{\widetilde{\Psi}(\widetilde{\Phi}^i(x))}$ then for a $1$-Lipschitz function $f$  we obtain
\begin{eqnarray}
  \vert\int\limits_{[0,1]}f\; d\widetilde{\omega}_{N,x}-\int\limits_{[0,1]}f\;d\omega_{N,x}\vert &\leq& \frac{1}{N}\sum_{i=0}^{N-1}\vert f(\widetilde{\Psi}(\widetilde{\Phi}^i(x)))-f(\Psi(\Phi^i(x)))\vert \nonumber\\
  &\leq & \frac{1}{N}\sum_{i=0}^{N-1}\vert \widetilde{\Phi}^{i+1}(x)-\widetilde{\Phi}^{i}(x)-\Phi^{i+1}(x)+\Phi^i(x)\vert<\frac{\varepsilon}{2} \nonumber
\end{eqnarray}
In other words, $d_{F}(\widetilde{\omega}_{N,x},\omega_{N,x})<\frac{\varepsilon}{2}$ for every $x\in [0,1]$, which means also that $d_{F}(\widetilde{\omega}_{N,x},\mu_{\Psi})<\varepsilon$.

Let now $x_0\in [0,1]$ be arbitrary. The orbit $\{\widetilde{\Phi}^{i}(x_0)\}_{i=0}^{\infty}$ can be divided into parts of length $N$: $\kappa_0=\{0,\widetilde{\Phi}(x_0),\dots, \widetilde{\Phi}^{N-1}(x_0)\}$, $\kappa_1=\{\widetilde{\Phi}^{N}(x_0),\widetilde{\Phi}^{N+1}(x_0),\dots, \widetilde{\Phi}^{2N-1}(x_0)\}$, $\kappa_2=\{\widetilde{\Phi}^{2N}(x_0),\widetilde{\Phi}^{2N+1}(x_0),\dots, \widetilde{\Phi}^{3N-1}(x_0)\}$ etc. and on each $\kappa_l$ part, $l\in\mathbb{N}\cup\{0\}$, the distribution $\widetilde{\omega}_{N,\widetilde{\Phi}^{Nl}(x_0)}$ approximates $\mu_{\Psi}$, i.e. $d_{F}(\widetilde{\omega}_{N,\widetilde{\Phi}^{Nl}(x_0)},\mu_{\Psi})<\varepsilon$. Define now the limit distribution: $\widetilde{\omega}_{x_0}=\lim_{n\to\infty}\frac{1}{n}\sum_{i=0}^{n-1}\delta_{\widetilde{\Psi}(\widetilde{\Phi}^i(x_0))}$. Note that the limit exists since if the rotation number $\widetilde{\varrho}$ of $\widetilde{\varphi}$ is irrational, then $\widetilde{\omega}_{x_0}$ equals the distribution $\widetilde{\mu}_{\widetilde{\Psi}}$ of the displacements of $\widetilde{\varphi}$ with respect to its invariant measure $\widetilde{\mu}$ and if $\widetilde{\varrho}\in\mathbb{Q}$, then $\widetilde{\omega}_{x_0}$ equals the displacement distribution on the periodic (finite) orbit to which the orbit of $x_0$ is attracted. Since on every $N$-part $\kappa_l$ the sample displacement distribution $\widetilde{\omega}_{N,\widetilde{\Phi}^{Nl}(x_0)}$ is $\varepsilon$-close in $d_F$-metric  to $\mu_{\psi}$, the limit distribution $\widetilde{\omega}_{x_0}$, which is obtained as averaging along the whole orbit of $x_0$, is also $\varepsilon$-close to $\mu_{\Psi}$:
\[d_F(\widetilde{\omega}_{x_0},\mu_{\Psi})<\varepsilon\]
Since $x_0$ was arbitrary, the proof is completed. \hbx

As the convergence under Fortet-Mourier metric implies weak convergence (cf. \cite{probmetric}) we conclude:
\begin{corollary}\label{wniosek100}
For every $x\in[0,1]$
\[\widetilde{\omega}_{x}\Longrightarrow \mu_{\Psi} \quad \textrm{as} \ \widetilde{\varphi}\to\varphi \ \textrm{in} \ C^0(S^1),\]
where $\widetilde{\omega}_{x}=\lim_{n\to\infty}\frac{1}{n}\sum_{i=0}^{n-1}\delta_{\widetilde{\Psi}(\widetilde{\Phi}^i(x))}$.
\end{corollary}

From the proofs of Theorem \ref{przyblizanie wymiernymi} and Corollary \ref{wniosek100} it follows that their statements are valid, naturally, also when one considers the invariant measure $\mu$ and the sample distribution of orbits $\varpi_{n,x}:=\frac{1}{n}\sum_{i=0}^{n-1}\delta_{\Phi^{i}(x)}$ instead of the distribution $\mu_{\Psi}$ and the sample displacement distribution $\omega_{n,x}$. Namely, we have that:
\begin{remark}
If $\varphi$ is a circle homeomorphism with irrational rotation number, then for every $\varepsilon>0$ there exists a neigbourhood $\mathcal{U}$ of $\varphi$ in $C^{0}(S^1)$ such that for every homeomorphism $\widetilde{\varphi}\in\mathcal{U}$ and $x\in[0,1]$ it holds that
\[d_{F}(\widetilde{\varpi}_x,\mu)<\varepsilon,\]
where $\widetilde{\varpi}_x:=\lim_{n\to\infty}\frac{1}{n}\sum_{i=0}^{n-1}\delta_{\widetilde{\Phi}^i(x)}$. Consequently,
\[\widetilde{\varpi}_{x}\Longrightarrow \mu \quad \textrm{as} \ \widetilde{\varphi}\to\varphi \ \textrm{in} \ C^0(S^1).\]
\end{remark}

\subsection{Regularity properties of the displacement sequence}
 In this section we will formulate some results concerning  regularity of the behaviour of the sequence $\{\eta_n(z)\}$ inside its interval/set of concentration. For this we have to introduce the following definition:
\begin{definition}\label{powracajacy ciag}
We say that a sequence $\{a_n\}$, $n\in\mathbb{N}$, is \emph{almost strongly recurrent} if it satisfies
$$\forall_{\varepsilon > 0} \ \exists_{N\in\mathbb{N}} \ \forall_{n\in\mathbb{N}} \ \forall_{k\in\mathbb{N}} \ \exists_{i\in\{0,1,...,N\}} \ |a_{n+k+i}-a_{n}|<\varepsilon
$$
\end{definition}
For an almost strongly recurrent sequence we require that for each $n$ the set of returns $\mathcal{R}:=\{r: \ r=k+i\}$ of $a_n$ to its $\varepsilon$-neighbourhood may depend on $n$, although for all $n$ its gaps are bounded by the same $N$. Note that in literature there is also a notion of an \emph{almost periodic sequence} (\cite{berg}), for which with given $k\in\mathbb{N}$, the index $i\in\{0,1,\dots, N\}$ can be chosen uniformly for all $n\in\mathbb{N}$ and thus almost periodicity of a sequence is a stronger property.

Let then $(X,\varphi)$ be a discrete dynamical system, where $(X,d)$ is a metric space and $\varphi\in C^0(X)$.
\begin{definition}\label{punkt prawieokresowy}{\bf (cf. \cite{gh3})}
A point $x\in X$ is \emph{almost periodic} if
$$
\forall_{\varepsilon>0} \ \exists_{N\in\mathbb{N}} \ \forall_{k\in\mathbb{N}} \ \exists_{i\in\{0,1,...,N\}} \ d(\varphi^{k+i}(x),x)<\varepsilon,
$$ i.e. when for any neighbourhood $U$ of $x$ the set of numbers $n$ such that $\varphi^n(x)\in U$ is relatively dense in $\mathbb{N}$.
\end{definition}
Sometimes an almost periodic point is called also  strongly recurrent. Recall the following theorem of W. Gottschalk:

\begingroup
\setcounter{tmp}{\value{theorem}}
\setcounter{theorem}{6} 
\renewcommand\thetheorem{\Alph{theorem}}
\begin{theorem}\label{twe}{\bf (cf. \cite{gh3})} Let $X$ be a compact metric space.

Then the closure of the orbit of any almost periodic point is a minimal set. Conversely, all points of any minimal set are almost periodic.
\end{theorem}

\endgroup

\setcounter{theorem}{\thetmp} 

\begin{proposition}\label{ciag strongly recurrent}
Let $\varphi:S^1\to S^1$ be a homeomorphisms with an irrational rotation number.

\begin{tabular}{p{0.8cm}p{13.3cm}}

  $\ref{ciag strongly recurrent}.1$ & If $\varphi$ is transitive, then for all $z\in S^1$ the displacement sequence $\{\eta_n(z)\}$ is almost strongly recurrent, i.e. \\
\empty & $\forall_{\varepsilon>0} \ \exists_{N\in\mathbb{N}} \ \forall_{n\in N} \ \forall_{k\in\mathbb{N} \cup \{0\}} \ \exists_{i\in\{0,1,..., N\}} \ |\eta_{n+k+i}(z)-\eta_{n}(z)|<\varepsilon$ \\
  $\ref{ciag strongly recurrent}.2$ & If $\varphi$ is not-transitive, then the sequence $\{\eta_n(z)\}$ is almost strongly recurrent for all $z\in \Delta$. \\
  \end{tabular}
\end{proposition}

The theorem will follow easily from the lemma
\begin{lemma}\label{lemat o prawie okresowym} Let $(X,d)$ be a compact metric space and $(X,\varphi)$, where $\varphi\in C^0(X)$, a discrete dynamical system with $X$ being a minimal invariant set.

 Given $x\in X$ and $\varepsilon>0$, the set of return-times of the orbit $\{\varphi^j(x)\}_{j=1}^{\infty}$ to the $\varepsilon$-neighbourhood of any point $\varphi^n(x)$ of the orbit of $x$ is relatively dense with gaps bounded uniformly for all $n$:
$$\forall_{\varepsilon>0} \ \exists_{N\in\mathbb{N}} \ \forall_{n\in \mathbb{N}\cup\{0\}} \ \forall_{k\in\mathbb{N}\cup\{0\}} \ \exists_{i\in \{1,2,...,N\}} \ d(\varphi^{n+k+i}(x),\varphi^n(x))<\varepsilon $$
\end{lemma}
\noindent{\bf Proof.} By compactness of $X$ there exists $m$ such that the open balls $B_l=B(\varphi^l(x),\varepsilon/2)$ for $l=0,1,...,m$ cover $X$. There exist numbers $N_0, N_1, N_2, ..., N_m$ such that for any $k\in\mathbb{N}\cup \{0\}$, if $\varphi^k(\varphi^l(x))\in B_l$ then also $\varphi^{k+i_l}(\varphi^l(x))\in B_l$ for some $i_l\in\{1,2,...,N_l\}$. In this case $d(\varphi^{l+k+i_l}(x),\varphi^{l}(x))<\varepsilon$. Setting $N=\max_{l\in\{0,1,...,m\}}N_l$ we obtain that
$$
\forall_{l\in\{0,1,...,m\}} \ \forall_{k\in\mathbb{N}\cup \{0\}} \ \exists_{i\in\{1,2,...,N\}} \ d(\varphi^{l+k+i}(x),\varphi^l(x))<\varepsilon
$$
Consider now an arbitrary point $\varphi^n(x)$, $n>m$, of the orbit $\{\varphi^j(x)\}_{j=1}^\infty$. There exists $l\in\{0,1,...,m\}$ such that $\varphi^n(x)\in B_l$ and $\varphi^n(x)=\varphi^j(\varphi^l(x))$ for some $j\in\mathbb{N}$, which is then a returning time of a point $\varphi^l(x)$ to $B_l$. Consequently, there exists $\widetilde{k}_1\in\{1,2,...,N\}$ such that $\varphi^{j+\widetilde{k}_1}(\varphi^{l}(x))=\varphi^{\widetilde{k}_1}(\varphi^n(x))\in B_l$. Similarly, there exists $\widetilde{k}_2\in\{1,2,...,N\}$ such that $\varphi^{j+\widetilde{k}_1+\widetilde{k}_2}(\varphi^{l}(x))=\varphi^{\widetilde{k}_1+\widetilde{k}_2}(\varphi^n(x))\in B_l$, etc. The sequence of numbers $k_r:=\widetilde{k}_1+\widetilde{k}_2+ ... +\widetilde{k}_r$, $r\in\mathbb{N}$, is relatively dense in $\mathbb{N}$ and
for every $r$ $d(\varphi^{k_r}(\varphi^n(x)),\varphi^n(x))<\varepsilon$. \hbx

\noindent{\bf Proof of Proposition \ref{ciag strongly recurrent}.} Let $z\in \Delta$ (where  $\Delta=S^1$ if $\varphi$ is transitive) and $\varepsilon>0$ be arbitrary.  There exists $\delta<\varepsilon/2$ such that for every $z_1, z_2\in S^1$ we have $\vert \varphi(z_1)-\varphi(z_2)\vert<\varepsilon$ whenever $\vert z_1-z_2\vert<\delta$. On the account of Lemma \ref{lemat o prawie okresowym} there exists $N$ such that for all integers $n>0$ and $k\geq 0$ and some $i\in\{1,2,...,N\}$ we have  $\vert\varphi^{(n-1)+k+i}(z)-\varphi^{n-1}(z)\vert<\delta$. Then
$$\vert \varphi^{n+k+i}(z)-\varphi^{n+k+i-1}(z)-\varphi^n(z)+\varphi^{n-1}(z)\vert<\varepsilon/2+\delta<\varepsilon,$$
which implies, with a little effort, that $\vert \eta_{n+k+i}(z)-\eta_n(z)\vert <\varepsilon$. \hbx

\begin{remark} Note that Proposition \ref{ciag strongly recurrent} applies in more general setting where $(X,\varphi)$ is a minimal dynamical system on a compact metric space $(X,d)$ and the displacement sequence is defined simply as $\eta_n(x)\colon = d(\varphi^n(x),\varphi^{n-1}(x))$.
\end{remark}

\section{Discussion}

In the end let us make a few comments.

\begin{remark} Notice that our results concerning the displacement sequence $\eta_n(z)$, $n\in \mathbb{N}$, are also valid when one considers $n\to -\infty$ (or simply $n\in \mathbb{Z}$).
\end{remark}

\begin{remark} The potential continuity of the mapping $\varphi\mapsto\gamma$ from the topology of $C^{>=1}(S^1)$ to $C^1(S^1)$-topology (if true, probably under some further constraints on $\varphi$ or $\varrho(\varphi)$) would mean that for such diffeomorphisms with irrational rotation numbers, being enough close in appropriate topology, the densities of the invariant measures are uniformly close. This could serve for deriving more refined conditions for convergence of the distributions $\mu^{(n)}_{\Psi_n}$ to $\mu_{\Psi}$.
\end{remark}

In the introduction we mentioned the connection between the displacement sequence of a circle homeomorphism and the interspike-intervals of periodically driven integrate-and-fire models. However, majority of our results can be deduced easier for these models, especially for the so-called Leaky Integrate-and-Fire or Perfect Integrator, where for the last one we have explicit analytical formulas of the conjugating homeomorphism of the firing map and the invariant measure (cf. \cite{brette1,wmjs1}). Nevertheless, our results, in particular Theorem \ref{przyblizanie wymiernymi} and Corollary \ref{wniosek100}, rigorously explain numerical results obtained by other authors concerning  distribution of interspike-intervals for periodically driven integrate-and-fire models (see for example ISI histograms in \cite{keener1}). We roughly formulate it as follows:
\begin{remark}
Consider an integrate-and-fire system $\dot{x}=f_{\lambda}(t,x)$, where $\lambda\in\Lambda\subset\mathbb{R}^k$ is a vector of parameters that smoothly parameterizes the family $\{f_\lambda\}_{\lambda\in\Lambda}$ of functions $1$-periodic in $t$  and which induces firing maps $\{\Phi_{\lambda}\}_{\lambda\in\Lambda}$ being lifts of circle homeomorphisms. If $\lambda_0$ is a parameter value for which the firing rate is irrational, then there exists a range of parameters $\Lambda_0\subset \Lambda$ such that $\lambda_0\in\Lambda_0$ and for every $\lambda\in\Lambda_0$ the sample interspike-intervals distributions are arbitrary close (in $d_F$-metric) to the  interspike-interval distribution of $\Phi_{\lambda_0}$.
\end{remark}
We remark that our result covers also non linear integrate-and-fire models.

\subsection*{Acknowledgements} The first author was supported  by national research grants NN 201 373236  and NCN 2011/03/B\\/ST1/04533 and the second author by National Science Centre grant DEC-2011/01/N/ST1/02698.

We are indebted to Micha\l $\ $Rams for pointing us an idea of formulating the result in the kind of Theorem \ref{przyblizanie wymiernymi} together with the sketch of the proof.


\medskip
WAC\L AW MARZANTOWICZ, Faculty of Mathematics and Computer Sci., \, Adam Mickiewicz University of Pozna\'n, ul. Umultowska 87,\, 61-614 Pozna{\'n},
Poland\\
\emph{E-mail address}: marzan@amu.edu.pl\\

\vspace{0.4cm}
\noindent JUSTYNA SIGNERSKA, Institute of Mathematics, Polish Academy of Sciences, ul.\'Sniadeckich 8, 00-956 Warszawa, Poland;\\
Faculty of Applied Physics and Mathematics, Gda\'nsk University of Technology, ul. Narutowicza 11/12,\, 80-233 Gda\'nsk, Poland\\
\emph{E-mail address}: j.signerska@impan.pl
\end{document}